\documentclass{gtart}

\def\ifplaintex{\expandafter\ifx\csname documentclass\endcsname\relax}

\ifplaintex 
\else
\fi

\expandafter\ifx\csname beginpicture\endcsname\relax
\expandafter\ifx\csname documentclass\endcsname\relax
\input pictex \else
\input prepictex \input pictex \input postpictex \fi\fi

\def\gt{{\mathsurround=0pt\it $\cal G\mskip-2mu$eometry \&\ 
$\cal T\!\!$opology}}        %  journal title in recommended style

\def\gtp{{\mathsurround=0pt\it $\cal G\mskip-2mu$eometry \&\ 
$\cal T\!\!$opology $\cal P\!$ublications}}  % GT publications

\def\lognumber#1{\def\thelognumber{#1}}
\def\volumenumber#1{\def\thevolumenumber{#1}}
\def\papernumber#1{\def\thepapernumber{#1}}
\def\volumeyear#1{\def\thevolumeyear{#1}}

\def\pagenumbers#1#2{\def\startpage{#1}\def\finishpage{#2}}
\def\published#1{\def\publishdate{#1}}
\def\proposed#1{\def\theproposer{#1}}
\def\seconded#1{\def\theseconders{#1}}
\def\received#1{\def\receiveddate{#1}}
\def\revised#1{\def\reviseddate{#1}}
\def\accepted#1{\def\accepteddate{#1}}
\def\asciititle#1{\def\theasciititle{#1}}

\long\def\asciiabstract#1{\long\def\theasciiabstract{#1}}

\let\\\par\let\thelognumber\relax
\let\thevolumenumber\relax\let\thepapernumber\relax
\let\thevolumeyear\relax\let\thesamplenumber\relax\let\startpage\relax
\let\finishpage\relax\let\publishdate\relax\let\receiveddate\relax
\let\reviseddate\relax\let\accepteddate\relax\let\theasciititle\relax
\let\theasciiauthors\relax
\let\theasciiabstract\relax
\let\theasciiemail\relax\let\theshortauthors\relax\let\theshorttitle\relax

\long\def\maketitlep{   % start of definition of \maketitlep

\count0=\startpage

\gt\hfill      %   Journal title (top left) 
\beginpicture
\setcoordinatesystem units <0.33truein, 0.33truein> point at 2.2 0.9
\setplotsymbol ({$\cal G$})
\plotsymbolspacing=9truept
\circulararc 315 degrees from 0 1 center at 0 0
\setplotsymbol ({$\cal T$})
\circulararc 315 degrees from 1 -1 center at 1 0
\endpicture
\break
{\small\ifx\thesamplenumber\relax % sample?  
Volume \else Sample
\fi\thevolumenumber\ (\thevolumeyear)
\startpage--\finishpage\nl
Published: \publishdate}
\vglue 0.5truein plus 0.4fil minus 0.1truein

{\parskip=0pt\leftskip 0pt plus 1fil\def\\{\par\smallskip}{\ifplaintex\large
\else\Large\fi\bf\thetitle}\par\medskip}   

\vglue 0pt plus 0.1fil 

{\parskip=0pt\leftskip 0pt plus 1fil\def\\{\par}{\sc\theauthors}
\par\medskip}

\vglue 0pt plus 0.1fil 

{\small\parskip=0pt\let\newline\\
{\leftskip 0pt plus 1fil\def\\{\par}{\sl\theaddress}\par}
\expandafter\ifx\theemail\relax    % email address?
\relax\else\vglue 5pt plus 0.02fil minus 2pt\def\\{\stdspace{\rm 
and}\stdspace} 
\cl{Email:\stdspace\tt\theemail}\fi
\ifx\theurl\relax                  % URL given?
\relax\else\vglue 5pt plus 0.02fil minus 2pt\def\\{\stdspace{\rm 
and}\stdspace}
\cl{URL:\stdspace\tt\theurl}\fi\par}

\vglue 7pt plus 0.3fil minus 3pt

{\bf Abstract}
\vglue 5pt plus 0.1fil minus 2pt

\theabstract

\vglue 7pt plus 0.3fil minus 3pt

{\bf AMS Classification numbers}\quad Primary:\quad \theprimaryclass

Secondary:\quad \thesecondaryclass

\vglue 5pt plus 0.3fil minus 2pt

{\bf Keywords}\quad \thekeywords

\vglue 10pt plus 0.5fil minus 5pt

{\small  Proposed: \theproposer\hfill Received: \receiveddate\nl
Seconded: \theseconders\hfill 
\ifx\reviseddate\relax                         % paper revised?
Accepted: \accepteddate                        % no
\else
Revised: \reviseddate                          % yes
\fi}
\eject
}       %  end of definition of \maketitlep

\let\maketitlepage\maketitlep
\let\maketitle\maketitlepage

\font\phead=cmsl9 scaled 950
\font=cmsl9 scaled 1050
\font\pnum=cmbx10 scaled 913
\font\lnum=cmbx10 
\font\pfoot=cmsl9 scaled 950
\font=cmsl9 scaled 1050
\ifplaintex
\headline{\vbox to 0pt{\vskip -4.5mm\line{\small\phead\ifnum
\count0=\startpage ISSN 1364-0380 (on line)
1465-3060 (printed) \hfill {\pnum\folio}\else\ifodd\count0\def\\{ }% 
\ifx\theshorttitle\relax\thetitle\else\theshorttitle\fi\hfill{\pnum\folio}
\else\def\\{ and }{\pnum\folio}\hfill\ifx\theshortauthors\relax\theauthors
\else\theshortauthors\fi\fi\fi}\vss}}
\footline{\vbox to 0pt{\vglue 0mm\line{\small\pfoot\ifnum\count0=\startpage
\copyright\ \gtp\hfill\else
\gt, Volume \thevolumenumber\ (\thevolumeyear)\hfill\fi}\vss
}}
\else
\makeatletter
\def\@oddhead{{\small\lhead\ifnum\count0=\startpage ISSN 1364-0380 (on line)
1465-3060 (printed) \hfill {\lnum\number\count0}\else\ifodd\count0
\def\\{ }\ifx\theshorttitle\relax \thetitle \else\theshorttitle\fi\hfill
{\lnum\number\count0}\else\def\\{ and }{\lnum\number\count0}
\hfill\ifx\theshortauthors\relax 
\theauthors\else\theshortauthors\fi\fi\fi}}\def\@evenhead{\@oddhead}
\def\@oddfoot{\small\lfoot\ifnum\count0=\startpage\copyright\ \gtp\hfill\else
\gt, Volume \thevolumenumber\ (\thevolumeyear)\hfill\fi}
\def\@evenfoot{\@oddfoot}
\makeatother
\fi

\newwrite\gtoutfile
\long\gdef\makeheadfile{  %%% start of definition of \makeheadfile
{\def\\{, }\def\s{ }
\immediate\openout\gtoutfile head.xxx
\immediate\write\gtoutfile{To: math@arxiv.org}
\immediate\write\gtoutfile{Subject: put or rep NNNNN:pppp}
\immediate\write\gtoutfile{--text follows this line--}
\immediate\write\gtoutfile{Proxy-for: \ifx\theasciiauthors\relax
\theauthors\else\theasciiauthors\fi\s<\ifx\theasciiemail\relax\theemail\else\theasciiemail\fi>}
\immediate\write\gtoutfile{\noexpand\\}
\immediate\write\gtoutfile{Authors: \ifx\theasciiauthors\relax
\theauthors\else\theasciiauthors\fi}
\immediate\write\gtoutfile{Title: \ifx\theasciititle\relax
\thetitle\else\theasciititle\fi}
\immediate\write\gtoutfile{Subj-class: GT or SG or MG etc}
\immediate\write\gtoutfile{MSC-class: \theprimaryclass\ifx\thesecondaryclass\relax\else, \thesecondaryclass\fi}
\immediate\write\gtoutfile{Journal-ref: Geom. Topol. \thevolumenumber
(\thevolumeyear) \startpage-\finishpage}
\immediate\write\gtoutfile{Comments: Published by Geometry and Topology at}
\immediate\write\gtoutfile{\s\s http://www.maths.warwick.ac.uk/gt/GTVol\thevolumenumber/paper\thepapernumber.abs.html}
\immediate\write\gtoutfile{\noexpand\\}
\immediate\write\gtoutfile{}
\ifx\theasciiabstract\relax
\immediate\write\gtoutfile{\theabstract}\else
\immediate\write\gtoutfile{\theasciiabstract}\fi
\immediate\write\gtoutfile{}
\immediate\write\gtoutfile{\noexpand\\}
\immediate\write\gtoutfile{}
\immediate\closeout\gtoutfile}}  %%% end of definition of \makeheadfile

\def\maketitlepage{\maketitlep\makeheadfile}
\let\maketitle\maketitlepage

\def\ifplaintex{\expandafter\ifx\csname documentclass\endcsname\relax}

\ifplaintex 
\else
\fi

\expandafter\ifx\csname beginpicture\endcsname\relax
\expandafter\ifx\csname documentclass\endcsname\relax
\input pictex \else
\input prepictex \input pictex \input postpictex \fi\fi

\def\gt{{\mathsurround=0pt\it $\cal G\mskip-2mu$eometry \&\ 
$\cal T\!\!$opology}}        %  journal title in recommended style

\def\gtp{{\mathsurround=0pt\it $\cal G\mskip-2mu$eometry \&\ 
$\cal T\!\!$opology $\cal P\!$ublications}}  % GT publications

\def\lognumber#1{\def\thelognumber{#1}}
\def\volumenumber#1{\def\thevolumenumber{#1}}
\def\papernumber#1{\def\thepapernumber{#1}}
\def\volumeyear#1{\def\thevolumeyear{#1}}

\def\pagenumbers#1#2{\def\startpage{#1}\def\finishpage{#2}}
\def\published#1{\def\publishdate{#1}}
\def\proposed#1{\def\theproposer{#1}}
\def\seconded#1{\def\theseconders{#1}}
\def\received#1{\def\receiveddate{#1}}
\def\revised#1{\def\reviseddate{#1}}
\def\accepted#1{\def\accepteddate{#1}}
\def\asciititle#1{\def\theasciititle{#1}}

\long\def\asciiabstract#1{\long\def\theasciiabstract{#1}}

\let\\\par\let\thelognumber\relax
\let\thevolumenumber\relax\let\thepapernumber\relax
\let\thevolumeyear\relax\let\thesamplenumber\relax\let\startpage\relax
\let\finishpage\relax\let\publishdate\relax\let\receiveddate\relax
\let\reviseddate\relax\let\accepteddate\relax\let\theasciititle\relax
\let\theasciiauthors\relax
\let\theasciiabstract\relax
\let\theasciiemail\relax\let\theshortauthors\relax\let\theshorttitle\relax

\long\def\maketitlep{   % start of definition of \maketitlep

\count0=\startpage

\gt\hfill      %   Journal title (top left) 
\beginpicture
\setcoordinatesystem units <0.33truein, 0.33truein> point at 2.2 0.9
\setplotsymbol ({$\cal G$})
\plotsymbolspacing=9truept
\circulararc 315 degrees from 0 1 center at 0 0
\setplotsymbol ({$\cal T$})
\circulararc 315 degrees from 1 -1 center at 1 0
\endpicture
\break
{\small\ifx\thesamplenumber\relax % sample?  
Volume \else Sample
\fi\thevolumenumber\ (\thevolumeyear)
\startpage--\finishpage\nl
Published: \publishdate}
\vglue 0.5truein plus 0.4fil minus 0.1truein

{\parskip=0pt\leftskip 0pt plus 1fil\def\\{\par\smallskip}{\ifplaintex\large
\else\Large\fi\bf\thetitle}\par\medskip}   

\vglue 0pt plus 0.1fil 

{\parskip=0pt\leftskip 0pt plus 1fil\def\\{\par}{\sc\theauthors}
\par\medskip}

\vglue 0pt plus 0.1fil 

{\small\parskip=0pt\let\newline\\
{\leftskip 0pt plus 1fil\def\\{\par}{\sl\theaddress}\par}
\expandafter\ifx\theemail\relax    % email address?
\relax\else\vglue 5pt plus 0.02fil minus 2pt\def\\{\stdspace{\rm 
and}\stdspace} 
\cl{Email:\stdspace\tt\theemail}\fi
\ifx\theurl\relax                  % URL given?
\relax\else\vglue 5pt plus 0.02fil minus 2pt\def\\{\stdspace{\rm 
and}\stdspace}
\cl{URL:\stdspace\tt\theurl}\fi\par}

\vglue 7pt plus 0.3fil minus 3pt

{\bf Abstract}
\vglue 5pt plus 0.1fil minus 2pt

\theabstract

\vglue 7pt plus 0.3fil minus 3pt

{\bf AMS Classification numbers}\quad Primary:\quad \theprimaryclass

Secondary:\quad \thesecondaryclass

\vglue 5pt plus 0.3fil minus 2pt

{\bf Keywords}\quad \thekeywords

\vglue 10pt plus 0.5fil minus 5pt

{\small  Proposed: \theproposer\hfill Received: \receiveddate\nl
Seconded: \theseconders\hfill 
\ifx\reviseddate\relax                         % paper revised?
Accepted: \accepteddate                        % no
\else
Revised: \reviseddate                          % yes
\fi}
\eject
}       %  end of definition of \maketitlep

\let\maketitlepage\maketitlep
\let\maketitle\maketitlepage

\font\phead=cmsl9 scaled 950
\font=cmsl9 scaled 1050
\font\pnum=cmbx10 scaled 913
\font\lnum=cmbx10 
\font\pfoot=cmsl9 scaled 950
\font=cmsl9 scaled 1050
\ifplaintex
\headline{\vbox to 0pt{\vskip -4.5mm\line{\small\phead\ifnum
\count0=\startpage ISSN 1364-0380 (on line)
1465-3060 (printed) \hfill {\pnum\folio}\else\ifodd\count0\def\\{ }% 
\ifx\theshorttitle\relax\thetitle\else\theshorttitle\fi\hfill{\pnum\folio}
\else\def\\{ and }{\pnum\folio}\hfill\ifx\theshortauthors\relax\theauthors
\else\theshortauthors\fi\fi\fi}\vss}}
\footline{\vbox to 0pt{\vglue 0mm\line{\small\pfoot\ifnum\count0=\startpage
\copyright\ \gtp\hfill\else
\gt, Volume \thevolumenumber\ (\thevolumeyear)\hfill\fi}\vss
}}
\else
\makeatletter
\def\@oddhead{{\small\lhead\ifnum\count0=\startpage ISSN 1364-0380 (on line)
1465-3060 (printed) \hfill {\lnum\number\count0}\else\ifodd\count0
\def\\{ }\ifx\theshorttitle\relax \thetitle \else\theshorttitle\fi\hfill
{\lnum\number\count0}\else\def\\{ and }{\lnum\number\count0}
\hfill\ifx\theshortauthors\relax 
\theauthors\else\theshortauthors\fi\fi\fi}}\def\@evenhead{\@oddhead}
\def\@oddfoot{\small\lfoot\ifnum\count0=\startpage\copyright\ \gtp\hfill\else
\gt, Volume \thevolumenumber\ (\thevolumeyear)\hfill\fi}
\def\@evenfoot{\@oddfoot}
\makeatother
\fi

\newwrite\gtoutfile
\long\gdef\makeheadfile{  %%% start of definition of \makeheadfile
{\def\\{, }\def\s{ }
\immediate\openout\gtoutfile head.xxx
\immediate\write\gtoutfile{To: math@arxiv.org}
\immediate\write\gtoutfile{Subject: put or rep NNNNN:pppp}
\immediate\write\gtoutfile{--text follows this line--}
\immediate\write\gtoutfile{Proxy-for: \ifx\theasciiauthors\relax
\theauthors\else\theasciiauthors\fi\s<\ifx\theasciiemail\relax\theemail\else\theasciiemail\fi>}
\immediate\write\gtoutfile{\noexpand\\}
\immediate\write\gtoutfile{Authors: \ifx\theasciiauthors\relax
\theauthors\else\theasciiauthors\fi}
\immediate\write\gtoutfile{Title: \ifx\theasciititle\relax
\thetitle\else\theasciititle\fi}
\immediate\write\gtoutfile{Subj-class: GT or SG or MG etc}
\immediate\write\gtoutfile{MSC-class: \theprimaryclass\ifx\thesecondaryclass\relax\else, \thesecondaryclass\fi}
\immediate\write\gtoutfile{Journal-ref: Geom. Topol. \thevolumenumber
(\thevolumeyear) \startpage-\finishpage}
\immediate\write\gtoutfile{Comments: Published by Geometry and Topology at}
\immediate\write\gtoutfile{\s\s http://www.maths.warwick.ac.uk/gt/GTVol\thevolumenumber/paper\thepapernumber.abs.html}
\immediate\write\gtoutfile{\noexpand\\}
\immediate\write\gtoutfile{}
\ifx\theasciiabstract\relax
\immediate\write\gtoutfile{\theabstract}\else
\immediate\write\gtoutfile{\theasciiabstract}\fi
\immediate\write\gtoutfile{}
\immediate\write\gtoutfile{\noexpand\\}
\immediate\write\gtoutfile{}
\immediate\closeout\gtoutfile}}  %%% end of definition of \makeheadfile

\def\maketitlepage{\maketitlep\makeheadfile}
\let\maketitle\maketitlepage

\def\ifplaintex{\expandafter\ifx\csname documentclass\endcsname\relax}

\ifplaintex 
\else
\fi

\expandafter\ifx\csname beginpicture\endcsname\relax
\expandafter\ifx\csname documentclass\endcsname\relax
\input pictex \else
\input prepictex \input pictex \input postpictex \fi\fi

\def\gt{{\mathsurround=0pt\it $\cal G\mskip-2mu$eometry \&\ 
$\cal T\!\!$opology}}        %  journal title in recommended style

\def\gtp{{\mathsurround=0pt\it $\cal G\mskip-2mu$eometry \&\ 
$\cal T\!\!$opology $\cal P\!$ublications}}  % GT publications

\def\lognumber#1{\def\thelognumber{#1}}
\def\volumenumber#1{\def\thevolumenumber{#1}}
\def\papernumber#1{\def\thepapernumber{#1}}
\def\volumeyear#1{\def\thevolumeyear{#1}}

\def\pagenumbers#1#2{\def\startpage{#1}\def\finishpage{#2}}
\def\published#1{\def\publishdate{#1}}
\def\proposed#1{\def\theproposer{#1}}
\def\seconded#1{\def\theseconders{#1}}
\def\received#1{\def\receiveddate{#1}}
\def\revised#1{\def\reviseddate{#1}}
\def\accepted#1{\def\accepteddate{#1}}
\def\asciititle#1{\def\theasciititle{#1}}

\long\def\asciiabstract#1{\long\def\theasciiabstract{#1}}

\let\\\par\let\thelognumber\relax
\let\thevolumenumber\relax\let\thepapernumber\relax
\let\thevolumeyear\relax\let\thesamplenumber\relax\let\startpage\relax
\let\finishpage\relax\let\publishdate\relax\let\receiveddate\relax
\let\reviseddate\relax\let\accepteddate\relax\let\theasciititle\relax
\let\theasciiauthors\relax
\let\theasciiabstract\relax
\let\theasciiemail\relax\let\theshortauthors\relax\let\theshorttitle\relax

\long\def\maketitlep{   % start of definition of \maketitlep

\count0=\startpage

\gt\hfill      %   Journal title (top left) 
\beginpicture
\setcoordinatesystem units <0.33truein, 0.33truein> point at 2.2 0.9
\setplotsymbol ({$\cal G$})
\plotsymbolspacing=9truept
\circulararc 315 degrees from 0 1 center at 0 0
\setplotsymbol ({$\cal T$})
\circulararc 315 degrees from 1 -1 center at 1 0
\endpicture
\break
{\small\ifx\thesamplenumber\relax % sample?  
Volume \else Sample
\fi\thevolumenumber\ (\thevolumeyear)
\startpage--\finishpage\nl
Published: \publishdate}
\vglue 0.5truein plus 0.4fil minus 0.1truein

{\parskip=0pt\leftskip 0pt plus 1fil\def\\{\par\smallskip}{\ifplaintex\large
\else\Large\fi\bf\thetitle}\par\medskip}   

\vglue 0pt plus 0.1fil 

{\parskip=0pt\leftskip 0pt plus 1fil\def\\{\par}{\sc\theauthors}
\par\medskip}

\vglue 0pt plus 0.1fil 

{\small\parskip=0pt\let\newline\\
{\leftskip 0pt plus 1fil\def\\{\par}{\sl\theaddress}\par}
\expandafter\ifx\theemail\relax    % email address?
\relax\else\vglue 5pt plus 0.02fil minus 2pt\def\\{\stdspace{\rm 
and}\stdspace} 
\cl{Email:\stdspace\tt\theemail}\fi
\ifx\theurl\relax                  % URL given?
\relax\else\vglue 5pt plus 0.02fil minus 2pt\def\\{\stdspace{\rm 
and}\stdspace}
\cl{URL:\stdspace\tt\theurl}\fi\par}

\vglue 7pt plus 0.3fil minus 3pt

{\bf Abstract}
\vglue 5pt plus 0.1fil minus 2pt

\theabstract

\vglue 7pt plus 0.3fil minus 3pt

{\bf AMS Classification numbers}\quad Primary:\quad \theprimaryclass

Secondary:\quad \thesecondaryclass

\vglue 5pt plus 0.3fil minus 2pt

{\bf Keywords}\quad \thekeywords

\vglue 10pt plus 0.5fil minus 5pt

{\small  Proposed: \theproposer\hfill Received: \receiveddate\nl
Seconded: \theseconders\hfill 
\ifx\reviseddate\relax                         % paper revised?
Accepted: \accepteddate                        % no
\else
Revised: \reviseddate                          % yes
\fi}
\eject
}       %  end of definition of \maketitlep

\let\maketitlepage\maketitlep
\let\maketitle\maketitlepage

\font\phead=cmsl9 scaled 950
\font=cmsl9 scaled 1050
\font\pnum=cmbx10 scaled 913
\font\lnum=cmbx10 
\font\pfoot=cmsl9 scaled 950
\font=cmsl9 scaled 1050
\ifplaintex
\headline{\vbox to 0pt{\vskip -4.5mm\line{\small\phead\ifnum
\count0=\startpage ISSN 1364-0380 (on line)
1465-3060 (printed) \hfill {\pnum\folio}\else\ifodd\count0\def\\{ }% 
\ifx\theshorttitle\relax\thetitle\else\theshorttitle\fi\hfill{\pnum\folio}
\else\def\\{ and }{\pnum\folio}\hfill\ifx\theshortauthors\relax\theauthors
\else\theshortauthors\fi\fi\fi}\vss}}
\footline{\vbox to 0pt{\vglue 0mm\line{\small\pfoot\ifnum\count0=\startpage
\copyright\ \gtp\hfill\else
\gt, Volume \thevolumenumber\ (\thevolumeyear)\hfill\fi}\vss
}}
\else
\makeatletter
\def\@oddhead{{\small\lhead\ifnum\count0=\startpage ISSN 1364-0380 (on line)
1465-3060 (printed) \hfill {\lnum\number\count0}\else\ifodd\count0
\def\\{ }\ifx\theshorttitle\relax \thetitle \else\theshorttitle\fi\hfill
{\lnum\number\count0}\else\def\\{ and }{\lnum\number\count0}
\hfill\ifx\theshortauthors\relax 
\theauthors\else\theshortauthors\fi\fi\fi}}\def\@evenhead{\@oddhead}
\def\@oddfoot{\small\lfoot\ifnum\count0=\startpage\copyright\ \gtp\hfill\else
\gt, Volume \thevolumenumber\ (\thevolumeyear)\hfill\fi}
\def\@evenfoot{\@oddfoot}
\makeatother
\fi

\newwrite\gtoutfile
\long\gdef\makeheadfile{  %%% start of definition of \makeheadfile
{\def\\{, }\def\s{ }
\immediate\openout\gtoutfile head.xxx
\immediate\write\gtoutfile{To: math@arxiv.org}
\immediate\write\gtoutfile{Subject: put or rep NNNNN:pppp}
\immediate\write\gtoutfile{--text follows this line--}
\immediate\write\gtoutfile{Proxy-for: \ifx\theasciiauthors\relax
\theauthors\else\theasciiauthors\fi\s<\ifx\theasciiemail\relax\theemail\else\theasciiemail\fi>}
\immediate\write\gtoutfile{\noexpand\\}
\immediate\write\gtoutfile{Authors: \ifx\theasciiauthors\relax
\theauthors\else\theasciiauthors\fi}
\immediate\write\gtoutfile{Title: \ifx\theasciititle\relax
\thetitle\else\theasciititle\fi}
\immediate\write\gtoutfile{Subj-class: GT or SG or MG etc}
\immediate\write\gtoutfile{MSC-class: \theprimaryclass\ifx\thesecondaryclass\relax\else, \thesecondaryclass\fi}
\immediate\write\gtoutfile{Journal-ref: Geom. Topol. \thevolumenumber
(\thevolumeyear) \startpage-\finishpage}
\immediate\write\gtoutfile{Comments: Published by Geometry and Topology at}
\immediate\write\gtoutfile{\s\s http://www.maths.warwick.ac.uk/gt/GTVol\thevolumenumber/paper\thepapernumber.abs.html}
\immediate\write\gtoutfile{\noexpand\\}
\immediate\write\gtoutfile{}
\ifx\theasciiabstract\relax
\immediate\write\gtoutfile{\theabstract}\else
\immediate\write\gtoutfile{\theasciiabstract}\fi
\immediate\write\gtoutfile{}
\immediate\write\gtoutfile{\noexpand\\}
\immediate\write\gtoutfile{}
\immediate\closeout\gtoutfile}}  %%% end of definition of \makeheadfile

\def\maketitlepage{\maketitlep\makeheadfile}
\let\maketitle\maketitlepage

\def\ifplaintex{\expandafter\ifx\csname documentclass\endcsname\relax}

\ifplaintex 
\else
\fi

\expandafter\ifx\csname beginpicture\endcsname\relax
\expandafter\ifx\csname documentclass\endcsname\relax
\input pictex \else
\input prepictex \input pictex \input postpictex \fi\fi

\def\gt{{\mathsurround=0pt\it $\cal G\mskip-2mu$eometry \&\ 
$\cal T\!\!$opology}}        %  journal title in recommended style

\def\gtp{{\mathsurround=0pt\it $\cal G\mskip-2mu$eometry \&\ 
$\cal T\!\!$opology $\cal P\!$ublications}}  % GT publications

\def\lognumber#1{\def\thelognumber{#1}}
\def\volumenumber#1{\def\thevolumenumber{#1}}
\def\papernumber#1{\def\thepapernumber{#1}}
\def\volumeyear#1{\def\thevolumeyear{#1}}

\def\pagenumbers#1#2{\def\startpage{#1}\def\finishpage{#2}}
\def\published#1{\def\publishdate{#1}}
\def\proposed#1{\def\theproposer{#1}}
\def\seconded#1{\def\theseconders{#1}}
\def\received#1{\def\receiveddate{#1}}
\def\revised#1{\def\reviseddate{#1}}
\def\accepted#1{\def\accepteddate{#1}}
\def\asciititle#1{\def\theasciititle{#1}}

\long\def\asciiabstract#1{\long\def\theasciiabstract{#1}}

\let\\\par\let\thelognumber\relax
\let\thevolumenumber\relax\let\thepapernumber\relax
\let\thevolumeyear\relax\let\thesamplenumber\relax\let\startpage\relax
\let\finishpage\relax\let\publishdate\relax\let\receiveddate\relax
\let\reviseddate\relax\let\accepteddate\relax\let\theasciititle\relax
\let\theasciiauthors\relax
\let\theasciiabstract\relax
\let\theasciiemail\relax\let\theshortauthors\relax\let\theshorttitle\relax

\long\def\maketitlep{   % start of definition of \maketitlep

\count0=\startpage

\gt\hfill      %   Journal title (top left) 
\beginpicture
\setcoordinatesystem units <0.33truein, 0.33truein> point at 2.2 0.9
\setplotsymbol ({$\cal G$})
\plotsymbolspacing=9truept
\circulararc 315 degrees from 0 1 center at 0 0
\setplotsymbol ({$\cal T$})
\circulararc 315 degrees from 1 -1 center at 1 0
\endpicture
\break
{\small\ifx\thesamplenumber\relax % sample?  
Volume \else Sample
\fi\thevolumenumber\ (\thevolumeyear)
\startpage--\finishpage\nl
Published: \publishdate}
\vglue 0.5truein plus 0.4fil minus 0.1truein

{\parskip=0pt\leftskip 0pt plus 1fil\def\\{\par\smallskip}{\ifplaintex\large
\else\Large\fi\bf\thetitle}\par\medskip}   

\vglue 0pt plus 0.1fil 

{\parskip=0pt\leftskip 0pt plus 1fil\def\\{\par}{\sc\theauthors}
\par\medskip}

\vglue 0pt plus 0.1fil 

{\small\parskip=0pt\let\newline\\
{\leftskip 0pt plus 1fil\def\\{\par}{\sl\theaddress}\par}
\expandafter\ifx\theemail\relax    % email address?
\relax\else\vglue 5pt plus 0.02fil minus 2pt\def\\{\stdspace{\rm 
and}\stdspace} 
\cl{Email:\stdspace\tt\theemail}\fi
\ifx\theurl\relax                  % URL given?
\relax\else\vglue 5pt plus 0.02fil minus 2pt\def\\{\stdspace{\rm 
and}\stdspace}
\cl{URL:\stdspace\tt\theurl}\fi\par}

\vglue 7pt plus 0.3fil minus 3pt

{\bf Abstract}
\vglue 5pt plus 0.1fil minus 2pt

\theabstract

\vglue 7pt plus 0.3fil minus 3pt

{\bf AMS Classification numbers}\quad Primary:\quad \theprimaryclass

Secondary:\quad \thesecondaryclass

\vglue 5pt plus 0.3fil minus 2pt

{\bf Keywords}\quad \thekeywords

\vglue 10pt plus 0.5fil minus 5pt

{\small  Proposed: \theproposer\hfill Received: \receiveddate\nl
Seconded: \theseconders\hfill 
\ifx\reviseddate\relax                         % paper revised?
Accepted: \accepteddate                        % no
\else
Revised: \reviseddate                          % yes
\fi}
\eject
}       %  end of definition of \maketitlep

\let\maketitlepage\maketitlep
\let\maketitle\maketitlepage

\font\phead=cmsl9 scaled 950
\font=cmsl9 scaled 1050
\font\pnum=cmbx10 scaled 913
\font\lnum=cmbx10 
\font\pfoot=cmsl9 scaled 950
\font=cmsl9 scaled 1050
\ifplaintex
\headline{\vbox to 0pt{\vskip -4.5mm\line{\small\phead\ifnum
\count0=\startpage ISSN 1364-0380 (on line)
1465-3060 (printed) \hfill {\pnum\folio}\else\ifodd\count0\def\\{ }% 
\ifx\theshorttitle\relax\thetitle\else\theshorttitle\fi\hfill{\pnum\folio}
\else\def\\{ and }{\pnum\folio}\hfill\ifx\theshortauthors\relax\theauthors
\else\theshortauthors\fi\fi\fi}\vss}}
\footline{\vbox to 0pt{\vglue 0mm\line{\small\pfoot\ifnum\count0=\startpage
\copyright\ \gtp\hfill\else
\gt, Volume \thevolumenumber\ (\thevolumeyear)\hfill\fi}\vss
}}
\else
\makeatletter
\def\@oddhead{{\small\lhead\ifnum\count0=\startpage ISSN 1364-0380 (on line)
1465-3060 (printed) \hfill {\lnum\number\count0}\else\ifodd\count0
\def\\{ }\ifx\theshorttitle\relax \thetitle \else\theshorttitle\fi\hfill
{\lnum\number\count0}\else\def\\{ and }{\lnum\number\count0}
\hfill\ifx\theshortauthors\relax 
\theauthors\else\theshortauthors\fi\fi\fi}}\def\@evenhead{\@oddhead}
\def\@oddfoot{\small\lfoot\ifnum\count0=\startpage\copyright\ \gtp\hfill\else
\gt, Volume \thevolumenumber\ (\thevolumeyear)\hfill\fi}
\def\@evenfoot{\@oddfoot}
\makeatother
\fi

\newwrite\gtoutfile
\long\gdef\makeheadfile{  %%% start of definition of \makeheadfile
{\def\\{, }\def\s{ }
\immediate\openout\gtoutfile head.xxx
\immediate\write\gtoutfile{To: math@arxiv.org}
\immediate\write\gtoutfile{Subject: put or rep NNNNN:pppp}
\immediate\write\gtoutfile{--text follows this line--}
\immediate\write\gtoutfile{Proxy-for: \ifx\theasciiauthors\relax
\theauthors\else\theasciiauthors\fi\s<\ifx\theasciiemail\relax\theemail\else\theasciiemail\fi>}
\immediate\write\gtoutfile{\noexpand\\}
\immediate\write\gtoutfile{Authors: \ifx\theasciiauthors\relax
\theauthors\else\theasciiauthors\fi}
\immediate\write\gtoutfile{Title: \ifx\theasciititle\relax
\thetitle\else\theasciititle\fi}
\immediate\write\gtoutfile{Subj-class: GT or SG or MG etc}
\immediate\write\gtoutfile{MSC-class: \theprimaryclass\ifx\thesecondaryclass\relax\else, \thesecondaryclass\fi}
\immediate\write\gtoutfile{Journal-ref: Geom. Topol. \thevolumenumber
(\thevolumeyear) \startpage-\finishpage}
\immediate\write\gtoutfile{Comments: Published by Geometry and Topology at}
\immediate\write\gtoutfile{\s\s http://www.maths.warwick.ac.uk/gt/GTVol\thevolumenumber/paper\thepapernumber.abs.html}
\immediate\write\gtoutfile{\noexpand\\}
\immediate\write\gtoutfile{}
\ifx\theasciiabstract\relax
\immediate\write\gtoutfile{\theabstract}\else
\immediate\write\gtoutfile{\theasciiabstract}\fi
\immediate\write\gtoutfile{}
\immediate\write\gtoutfile{\noexpand\\}
\immediate\write\gtoutfile{}
\immediate\closeout\gtoutfile}}  %%% end of definition of \makeheadfile

\def\maketitlepage{\maketitlep\makeheadfile}
\let\maketitle\maketitlepage

\lognumber{167}
\volumenumber{5}\papernumber{24}\volumeyear{2001}
\pagenumbers{761}{797}
\received{23 February 2001}
\accepted{5 October 2001}
\proposed{Simon Donaldson}
\seconded{John Morgan, Tomasz Mrowka}
\revised{25 October 2001}
\published{25 October 2001}

\let\endpf\endproof
\def\S{Section }

\theoremstyle{remark}
\newtheorem{defn}{Definition}[section]
\newtheorem{remark}[defn]{Remark}

\newtheorem{definition}[defn]{Definition}

\theoremstyle{plain}
\newtheorem{lemma}[defn]{Lemma}

\newtheorem{theorem}[defn]{Theorem}
\newtheorem{proposition}[defn]{Proposition}
\newtheorem{corollary}[defn]{Corollary}

\newtheorem{claim}[defn]{Claim}
\newtheorem{maintheorem}{Theorem}

\newcommand{\kahler}{K\"{a}hler\ }
\newcommand{\tr}{{\rm Tr\,}}

\newcommand{\cc}{{\bf C}}
\newcommand{\zz}{{\bf Z}}
\newcommand{\rr}{{\bf R}}
\newcommand{\nn}{{\bf N}}

\newcommand{\calr}{{\cal R}}
\newcommand{\cale}{{\cal E}}
\newcommand{\calc}{{\cal C}}

\newcommand{\caln}{{\cal N}}
\newcommand{\cala}{{\cal A}}
\newcommand{\calw}{{\cal W}}
\newcommand{\call}{{\cal L}}
\newcommand{\calh}{{\cal H}}

\newcommand{\calo}{{\cal O}}
\newcommand{\calz}{{\cal Z}}
\newcommand{\cald}{{\cal D}}
\newcommand{\calk}{{\cal K}}
\newcommand{\calg}{{\cal G}}
\newcommand{\calu}{{\cal U}}
\newcommand{\calv}{{\cal V}}
\newcommand{\rplus}{{\cal R}_+}
\newcommand{\rminus}{{\cal R}_-}
\newcommand{\cpone}{{\bf CP^1}}

\newcommand{\calm}{{\cal M}}
\newcommand{\cali}{{\cal I}}
\newcommand{\cyl}{\Sigma\times S^1\times\rr}
\newcommand{\dbar}{\bar{\partial}}
\newcommand{\poincare}{Poincar\'{e} }

\newcommand{\cs}{{\bf cs}}

\begin{document}

\title{Instantons on cylindrical manifolds\\and stable bundles}
\asciititle{Instantons on cylindrical manifolds and stable bundles}
\author{Brendan Owens}
\address{Department of Mathematics and Statistics\\McMaster 
University\\Hamilton, Ontario, Canada}
\email{owensb@icarus.math.mcmaster.ca}

\begin{abstract}
Let $\Sigma$ be a smooth complex curve, and let $S$ be the product ruled
surface $\Sigma \times \cpone$.
We prove a correspondence conjectured by Donaldson between
finite energy $U(2)$--instantons over $\Sigma\times S^1\times\rr$, and rank 2
holomorphic bundles over $S$ whose restrictions to 
$\Sigma\times\{0\},\Sigma\times\{\infty\}$ are stable.
\end{abstract}

\asciiabstract{
Let Sigma be a smooth complex curve, and let S be the product ruled
surface Sigma \times CP^1.
We prove a correspondence conjectured by Donaldson between finite 
energy U(2)-instantons over the cylinder Sigma \times S^1 \times R, 
and rank 2 holomorphic bundles over S whose restrictions to the 
divisors at infinity are stable.
}

\maketitlepage

\section{Introduction}
Let $\Sigma$ be a smooth complex curve and let $Y=\Sigma\times S^1$.  Let
$E\to\cyl$ be a rank 2 complex vector bundle pulled back from $\Sigma$, with
$c_1=1$ over $\Sigma$.

Denote by $\calm$ the moduli space of finite energy $U(2)$--instantons on $E$.
There is a natural involution on $\calm$ which we denote by $\cali$.  This is
essentially given by tensoring $E$ with the flat complex line bundle pulled
back from the nontrivial double cover of the $S^1$ factor.

The cylinder $\cyl$ may be compactified by adding a copy of $\Sigma$ at 
each end.  We denote these curves by $\Sigma_0$ and $\Sigma_\infty$.
The resulting compactification is the product ruled surface 
$S=\Sigma\times\cpone$.  Let $\calz$ denote the space of isomorphism classes
of rank 2 holomorphic bundles over $S$ with fixed degree 1 determinant over 
$\Sigma$ and whose restriction to $\Sigma_0$ and $\Sigma_\infty$ is stable.

The following result was conjectured by Donaldson in \cite{d2}:

\begin{maintheorem}\label{mainthm}
The space $\calm/\cali$ is naturally homeomorphic to $\calz$.
\end{maintheorem}

(The topologies on the spaces $\calm/\cali$ and $\calz$ are
quotients of the $C^\infty$ topologies on compact sets, see Sections
\ref{instantons} and \ref{stablebundles} for details.)

There are four main steps in the proof of Theorem \ref{mainthm}.  
First we define a map
$$\Psi\co \calm/\cali\to\calz.$$
Given $[A]\in \calm$ we obtain $\Psi([A])$
by extending the holomorphic bundle $(E,\dbar_A)$ over $\Sigma_0$ and
$\Sigma_\infty$.  The idea behind this extension is a simple one which
exploits the product structure of the 3--manifold $\Sigma\times S^1$.
Standard convergence results for finite energy instantons enable us to show
that the restriction of $(E,\dbar_A)$ to every fiber 
$\Sigma\times(\theta,t)\subset\Sigma\times S^1\times\rr$ is stable for
large $|t|$.  Thus the bundle on each end of the tube gives a holomorphic
map from a punctured disk into the moduli space $\caln(\Sigma)$ of stable
bundles on $\Sigma$; by
removable singularities we are able to extend this map over the puncture
and thus also extend the bundle over the divisors at infinity.

The remaining steps in the proof are to show that the map $\Psi$ is injective,
onto, and finally a homeomorphism.  The proof of injectivity follows an
argument of Donaldson based on a function comparing two Hermitian metrics.
We show
surjectivity by exhibiting a construction which gives a compatible 
$U(2)$--instanton on the restriction of a vector bundle $\cale$ to
$\Sigma\times S^1\times\rr\subset S$, where $\cale$ represents an element
of $\calz$.  This construction is based on Donaldson's theorem
that any stable bundle on a compact algebraic surface admits a 
Hermitian--Einstein connection, and on Uhlenbeck compactness for instantons.
We take a sequence of Hodge metrics on $S$ which converge to the cylinder
metric on
compact subsets of $\Sigma\times S^1\times\rr$.  Donaldson's theorem gives
a corresponding sequence of Hermitian--Einstein connections on $\cale$, which
has a weakly convergent subsequence yielding a $U(2)$--instanton $A$ on the
cylinder.  The stability of $\cale$ ensures that no energy is lost in this
process and we are able to show that $\Psi([A])=[\cale]$.

The last step in the proof is to show that $\Psi$ is a homeomorphism.
The key point here (also used in the proof of surjectivity)
is that for certain sequences of stable bundles on $S$, convergence on
a suitably chosen subset of $S$ actually implies global convergence.
\\

{\bf Acknowledgements}\qua This paper is a modified version of my Columbia 
University PhD thesis.  I am very grateful to my advisor, John Morgan, for
his invaluable help and advice.  
Thanks also to Michael Thaddeus, Robert Friedman,
Mehrzad Ajoodanian, Pedram Safari and Saso Strle for useful discussions.
I am grateful to the referee for pointing out an error in an earlier draft.

\section{Definitions and notation}
\label{defnot}
We let $Y$ denote the 3--manifold $\Sigma\times S^1$ as above.  Fix a \kahler
metric $g_\Sigma$ on $\Sigma$.  Let $\theta$ and $t$ be the standard
coordinates on $S^1$ and $\rr$ respectively.  Then
$$g_{cyl}=g_\Sigma+d\theta^2+dt^2$$
is a complete \kahler metric on $\cyl$.

The cylinder $\cyl$ embeds in the product ruled surface $S=\Sigma\times\cpone$
in such a way that the standard coordinate $z$ on $\cpone$ is given by
$$z=e^t e^{i\theta}.$$
Let $E_\Sigma$ denote the rank 2, degree 1 complex vector bundle over $\Sigma$,
with a chosen Hermitian metric.
Let $E_Y$ and $E$ denote the pullbacks of $E_\Sigma$ to $Y$ and 
$Y\times\rr$ respectively.

Given a unitary connection $A$ on $E$
we define the energy of $A$ to be
$$e(A)=\int_{Y \times \rr}|F_A^0|^2\,{\rm dvol}
=-\int_{Y \times \rr}\tr F_A^0\wedge*F_A^0$$
where $F_A^0$ is the trace-free part of the curvature.  This is equal to
the Yang--Mills energy of the $SO(3)$ connection obtained by projectivising
$A$ and $E$.

\begin{definition}
\label{instantondef}
A $U(2)$--instanton on $E$ is a
finite-energy projectively anti-self-dual connection on $E$, with fixed
central part.  That is, 
a unitary connection $A$, which with respect to the metric $g_{cyl}$, has
\begin{itemize}
\item trace-free part of curvature $F_A^0$ is anti-self-dual
\item finite energy
\item fixed central part.
\end{itemize}
For convenience, we choose the fixed central part to be pulled back from 
$\Sigma$, with $\tr F_A$ harmonic.
\end{definition}

We denote by $\calm$ the space of $U(2)$--instantons modulo gauge 
transformations which fix the central part.  
We denote by $[A]$ the element of $\calm$ 
represented by a $U(2)$--instanton $A$.
Let $\calm_e$ denote the set
of equivalence classes with energy a fixed real number $e$.  For $e<0$, $\calm_e$ is empty; for $e=0$, $\calm_e$
consists of equivalence classes of flat connections.  We will see in
Section \ref{instantons} that $e$ is always an integer multiple of 
$4\pi^2$.

The definition of $U(2)$--instanton
specifies that the central part of the connection is fixed and pulled back
from $\Sigma$.  This then determines a holomorphic structure on the degree 1
smooth complex line bundle on $\Sigma$; 
we denote this holomorphic line bundle by
$\call$.  We define $\cald$ to be the line bundle on $S$ determined by
$\call$ and a number $e\in4\pi^2\zz$ as follows:
\begin{equation}
\label{det}
\cald=
\left\{\begin{array}{ll}
\pi_1^*\call&\mbox{if}\quad e(A)\in8\pi^2\zz,\\
\pi_1^*\call\otimes\pi_2^*\calo_\cpone(1)&\mbox{otherwise}.
\end{array}\right.
\end{equation}
We also define an integer $c$, again determined by $\call,e$:
\begin{equation}
\label{c2}
c=\frac14c_1(\cald)^2+\frac{1}{8\pi^2}e(A).
\end{equation}

\begin{definition}
\label{zdef}
Fix a line bundle $\cald\to S$ and an integer $c$.
Let $\calz_{(\cald,c)}$ denote the set of isomorphism classes of rank 2
holomorphic bundles $\cale$ on $S$ satisfying
\begin{itemize}
\item ${\rm det}(\cale)\cong\cald$
\item $c_2(\cale)=c$
\item the restriction of $\cale$ to $\Sigma_0, \Sigma_\infty$ is stable.
\end{itemize}
\end{definition}

For a given bundle $\cale$ as above, let $[\cale]$ be the element of 
$\calz_{(\cald,c)}$ it represents.

\section{Finite energy instantons on the cylinder}
\label{instantons}
We start by studying the moduli space $\calm$ of $U(2)$--instantons on 
$Y\times\rr$ (defined in
Section \ref{defnot}).
The topology on this space is as follows:
a sequence $[A_n]$ in $\calm$ converges to a limit $[A_\infty]$, if and only
if for some representatives,
$A_n$ converges to $A_\infty$ in $C^\infty$ on compact subsets of $Y\times\rr$.

\begin{remark}Throughout this paper (unless explicitly stated otherwise)
unitary gauge equivalence refers only to gauge transformations which fix
the central part.  This is actually equivalent to considering the group of
even $SO(3)$ gauge transformations (see \cite{bd}).
\end{remark}

Projectively flat connections on $E_Y\to Y$ are unitary connections with a 
fixed central part whose curvature is central.
The moduli space of projectively flat connections
on $E_Y$ modulo gauge is denoted by $\calr$.  It consists of two
connected components $\rplus$ and $\rminus$, each of which is a copy of the
moduli space of projectively flat connections on $E_\Sigma$.  By the 
Narasimhan--Seshadri theorem this
moduli space is diffeomorphic to the moduli space 
$\caln(\Sigma)$ of rank 2, degree 1 stable
bundles over $\Sigma$ (with fixed determinant).  It is a smooth compact 
complex manifold of complex 
dimension $3g-3$.  The sets
$\rplus$ and $\rminus$ consist of those projectively flat connections whose
holonomy around the $S^1$ factor is $+{\bf 1}$ and $-{\bf 1}$ respectively.

We now recall some results about finite energy instantons
on cylindrical end manifolds.

Let $A$ be a $U(2)$--instanton on $E \rightarrow Y \times \rr$.  Fixing a
bundle isomorphism
$$\eta\co E_Y\times\rr\to E,$$
the pullback $\eta^*A$ gives a 1--parameter family of unitary connections
on $E_Y$, which we denote by $A_t$.
We will need to use the following results about the
behaviour of the gauge equivalence class $[A_t]$:

\begin{theorem}
\label{thmA}
The class $[A_t]$ converges to a projectively flat connection
$[A_{\pm\infty}]$ on $Y$ as $t \rightarrow \pm\infty$.
\end{theorem}

\begin{theorem}
\label{thmB}
For any gauge representatives $A_{\pm\infty}$ of the limits of $[A_t]$,
there are numbers $C,\delta>0$ and smooth gauge transformations $g_\pm$
such that the following estimates hold:
$$||g_\pm^*A-\pi_Y^*A_{\pm\infty}||_{L^2_2(Y\times [T-\frac12,T+\frac12])}
\le C e^{\mp\delta T}\quad\mbox{for}\,\pm T>>0.$$
\end{theorem}

The first of these results follows from Theorem 4.0.1 of \cite{mmr}, with 
$G=SO(3)$ (note that $\calm$ may equivalently be described as the space of
finite-energy $SO(3)$ instantons on the projectivisation of $E$, modulo
gauge transformations which lift to $SU(2)$).  The second follows from Theorem
5.2.2 of \cite{mmr} (see also Lemma 2.1.10 in \cite{mm}).

If $A$ is a $U(2)$--instanton converging to projectively flat limits
$A_{\pm\infty}$ as in Theorem \ref{thmA}, then the energy of $A$ is related to
the value of the Chern--Simons functional on the limits as follows:
$$\frac{1}{8\pi^2}e(A)=
\frac{1}{8\pi^2}\int_{Y \times \rr}\tr F_A^0\wedge F_A^0=
(\cs(A_{+\infty})-\cs(A_{-\infty}))\bmod\zz.$$
(In fact the second equality is the definition of the Chern--Simons functional,
 which is $\rr/\zz$--valued and is only defined up to an additive constant.)

Given a projectively flat connection $a$ on $\Sigma$,
the pullback $\pi_\Sigma^*a$ represents an element of $\rplus$.  
Trivialise the degree one $U(2)$ bundle over a disk
$D\subset \Sigma$ and over its complement, with transition function
$$
\left(\begin{array}{cc}
e^{i\phi} & 0\\
0 & 1 \end{array} \right)$$
along the boundary circle.
Then $g = \left(\begin{array}{cc}
e^{-i\theta} & 0\\
0 & 1 \end{array} \right)$ is a well-defined unitary gauge transformation
over $Y$.

Then the connection  
$$\cali_Y(\pi_\Sigma^*a)=
g^*(\pi_\Sigma^*a + \frac{i}{2} d\theta \cdot {\bf 1}_{E_Y})$$
represents an element of $\rminus$  (Adding the term
$\frac{i}{2} d\theta \cdot {\bf 1}_{E_Y}$ gives the desired holonomy around the
$S^1$ factor, the gauge transformation is required so that elements of $\rplus$
and $\rminus$ have the same central part.)  Note that $g$ is not in the
restricted gauge group of unitary gauge transformations which fix the central
part, but that the square of $\cali_Y$ is equivalent to applying the
gauge transformation $e^{i\theta}g^2$, which is.

The given map $\cali_Y$ is thus an involution on $\calr$ which 
interchanges the 
components $\rplus$ and $\rminus$.  It also changes the value of $\cs$ by 
$\frac12$ (see \cite{bd,ds}).
It follows that the energy $e(A)$ is an integer multiple of $4\pi^2$.
Also $\frac{1}{4\pi^2}e(A)$ is even if the limits $A_{\pm\infty}$ are in the
same component, odd otherwise.

Let $A\in\calm_e$ for some $e\in 4\pi^2\zz$.  Then $A$ defines a type $(0,1)$ 
operator $\dbar_A$ on $E$--valued forms which satisfies $\dbar_A{}^2=0$.  The 
Newlander--Nirenberg theorem implies that $\dbar_A$ defines a holomorphic
structure on the complex vector bundle $E$ (see eg \cite[page 46]{dk}).

We will now give the definition of an involution
$$\cali\co \calm_e\to\calm_e$$
and show that $(E,\dbar_A)\cong(E,\dbar_{\cali(A)})$ as holomorphic bundles.

\begin{lemma}
\label{involution}
There is an involution ${\cal I}\co \calm_e \to \calm_e$ which 
fixes the holomorphic structure on $E$ determined by a $U(2)$--instanton.
Its restriction to $Y$ is the involution $\cali_Y$ which switches
the components $\rplus,\rminus$ of the space $\calr$ of projectively flat
connections.
\end{lemma}

\proof
We define the map ${\cal I}$ by
$$A \mapsto g^*(A + \frac{i}{2} d\theta \cdot {\bf 1}_E),$$
where $g$ is a $U(2)$ gauge transformation to fix the central part as
required by the definition of $U(2)$--instanton.  (For example take
$g = \left(\begin{array}{cc}
e^{-i\theta} & 0\\
0 & 1 \end{array} \right)$ as above).

It is fairly clear that this defines an involution on $\calm$ whose
restriction to $Y$ is the involution $\cali_Y$.  It also clearly
preserves the energy.

It remains to see that it
preserves the holomorphic structure.  This follows from the 
$\bar{\partial}$--Poincar\'{e} lemma on $\cc - \{0\}$, 
applied to the $(0,1)$ part of
$\frac{i}{2} d\theta$.  This shows that the map
$$A \mapsto A + \frac{i}{2} d\theta \cdot {\bf 1}_E$$
is given by a (central) complex gauge transformation.

(In fact taking $h=|z|^{-\frac12}$ yields 
$h^{-1}\dbar h=\frac{i}{2} d\theta^{(0,1)}$, as the reader may verify.) \endpf

Using Theorem \ref{thmB} and Lemma \ref{involution} we can modify a 
$U(2)$--instanton in order to get good gauge representatives on each end of the
cylinder.

\begin{lemma}
\label{gaugereps}
Let $A$ be a $U(2)$--instanton on 
$E\to \Sigma \times S^1 \times \rr$ with energy \mbox{$e(A)=e\in4\pi^2\zz$}.  
Then there exists a projectively ASD connection ${\tilde A}$ on 
$E\to \Sigma \times S^1 \times \rr$ with the following properties:
\begin{itemize}
\item $(E,\dbar_A)\cong(E,\dbar_{{\tilde A}})$ as holomorphic bundles
\item $e(A)=e({\tilde A})$
\item the limits $[{\tilde A}_{\pm\infty}]$ on $Y$ are gauge equivalent (under the
full $U(2)$ gauge group) to the pullbacks
$\pi_\Sigma^*a_\pm$ of projectively flat connections $a_\pm$ on $\Sigma$.
\end{itemize}
\end{lemma}

\proof
There are 4 cases to consider, depending on whether the limits 
$[A_{\pm\infty}]$ are in $\rplus$ or $\rminus$.

The first two cases occur when both limits $[A_{\pm\infty}]$ are in the same
component of $\calr$; this occurs when $\frac{e}{4\pi^2}$ is even.

If both of $[A_{\pm\infty}]$ are in $\rplus$, then simply take ${\tilde A}=A$.

If both of $[A_{\pm\infty}]$ are in $\rminus$, then take ${\tilde A}=\cali(A)$.

The other two cases occur when $A$ has ``mixed limits'', ie, when
$\frac{e}{4\pi^2}$ is odd.

Suppose $[A_{-\infty}]\in\rplus$ and 
$[A_{+\infty}]\in\rminus$.  Then $[A_{-\infty}]$ is represented by the 
pullback 
$\pi_\Sigma^*a_{-\infty}$ of a projectively flat connection 
$a_{-\infty}$ on $\Sigma$.
The limit $[A_{+\infty}]$ is represented by 
$$A_{+\infty}=g^*(\pi_\Sigma^*a_{+\infty}  + \frac{i}{2} d\theta \cdot {\bf 1}_E)
=g^*(\pi_\Sigma^*a_{+\infty}) + \frac{i}{2} d\theta \cdot {\bf 1}_E,$$
with $g$ as in the proof of Lemma \ref{involution} and $a_{+\infty}$ a 
projectively flat connection on $\Sigma$.

Define $${\tilde A}=A-s(t)\frac{i}{2} d\theta \cdot {\bf 1}_E$$
where $s(t)$ is a smooth function on $\rr$ with
$$s(t)=
\left\{\begin{array}{ll}
0 &\mbox{if}\quad t<0,\\
1 &\mbox{if}\quad t>1.
\end{array}\right.$$
This ${\tilde A}$ satisfies the required conditions.  It is projectively ASD 
with the
same energy as $A$ since its trace-free part is the same as that of $A$.
The holomorphic bundle $(E,\dbar_{{\tilde A}})$ is equal to $(E,\dbar_A)$ 
tensored
by a line bundle pulled back from $\cc^*$; any such line bundle is 
holomorphically trivial.  By construction the limits 
$[{\tilde A}_{\pm\infty}]$ have 
the required property.

If $[A_{-\infty}]\in\rminus$ and $[A_{+\infty}]\in\rplus$, define
${\tilde A}={\widetilde {\cali(A)}}$.

Note that if $[A],[B]\in\calm_e$ then ${\tilde A}$ is gauge equivalent to
${\tilde B}$ by a unitary gauge transformation fixing the central part 
if and only if $[A]=[B]$ or $[A]=\cali([B])$.
\endpf

\begin{remark}
The preceding Lemma may be considered from the point of view of orbifold
bundles and connections.  A connection on $\cyl$ which is converging to a
flat connection with holonomy $-{\bf 1}$ around the $S^1$ factor will clearly
not extend over $\Sigma\times\cpone$.  We may think of such a connection as
a connection on the orbifold $\Sigma\times\cpone$ with $\zz_2$ acting on
$\cpone$ by $z\mapsto z^2$.  We can then produce a connection with $+{\bf 1}$
holonomy by tensoring with an orbifold line bundle connection with matching
holonomy limits.
\end{remark}

\section{Stable bundles on $\Sigma$}
In this section we briefly note two very important properties of the moduli
space $\caln(\Sigma)$ of rank 2, degree 1 
stable bundles on $\Sigma$ with fixed determinant.  For a detailed study of
this space see \cite{ab}, \cite{n}.

\label{finemoduli}
Firstly, $\caln(\Sigma)$ is a fine moduli space
(see \cite[\S9]{ab},\cite{n} for a proof, also \cite[\S4.2]{fm} for 
general discussion).
What this means is that there is a holomorphic bundle $\calu$ over 
$\Sigma \times \caln(\Sigma)$ with the following universal property.  
Let $M$ be any complex manifold.
For any holomorphic bundle $\calv$ over $\Sigma
\times M$, whose restriction to each $\Sigma \times \{{\rm point}\}$ 
is stable, there
exists a unique holomorphic map $f\co M \to \caln(\Sigma)$ and a 
line bundle $\calh$ on $M$
such that $\calv\cong({\rm Id}\times f)^*\calu \otimes \pi_2^*\calh$.

\label{slices}
Secondly, there are local slices for the action of the complex
gauge group on the space of stable $(0,1)$ connections 
on $E_\Sigma\to\Sigma$.  This means the following.  Let $\cala_0^{(0,1)}$
denote the space of all $(0,1)$ connections on $E_\Sigma$ with fixed
central part, and let $\calg_0^c=\Gamma(SL(E_\Sigma))$ denote the group
of complex gauge automorphisms fixing the central part.  
(In fact for the following discussion we take the $L^2_1$ completion of
$\cala_0^{(0,1)}$ and the $L^2_2$ completion of the complex gauge group.)
The group 
$\calg_0^c$ acts on $\cala_0^{(0,1)}$ with stabiliser $\pm{\bf 1}$.
The set of all holomorphic structures on $E_\Sigma$ is the quotient of
$\cala_0^{(0,1)}$ by this action.  This space is not Hausdorff.  If we let
$\cala^s\subset \cala_0^{(0,1)}$ denote the subset of stable holomorphic
structures, then the quotient of $\cala^s$ by $\calg_0^c$ is the moduli space
$\caln(\Sigma)$.

The differential
at any point $\dbar_0\in\cala_0^{(0,1)}$ of the complex gauge group action is 
given by the operator 
$$\dbar_0\co \Omega^0(\Sigma,{\rm End}\,E_\Sigma)\to
\Omega^{0,1}(\Sigma,{\rm End}\,E_\Sigma)$$
(restricted to trace-free forms).
The metrics on $\Sigma$ and $E_\Sigma$ enable us to define a formal adjoint 
$\dbar_0^*$.
Denote by
$$S_{\dbar_0}={\dbar_0}+{\rm Ker}\ \dbar_0^*$$ the affine slice in
 $\cala_0^{(0,1)}$.  This is an affine complex subspace which is
 transverse to the gauge orbit at the point $\dbar_0$.  Suppose that
 $\dbar_0$ gives a stable holomorphic structure.  The action of
 $\calg_0^c/\!_{\{\pm{\bf 1}\}}$ on $\cala^s$ is smooth, free and
 proper; it follows that there is a neighbourhood $U$ of $\dbar_0$ in
 $S_{\dbar_0}$, such that the natural map
$$U\times_{\pm{\bf 1}}\calg_0^c\to \cala_0^{(0,1)}$$
is a diffeomorphism onto its image.
(See \cite[pages 294--295 and page 300]{fm} 
for details.  The discussion there is for a base manifold of any dimension;
the dimension one case is simpler since every $(0,1)$ connection is 
integrable.)

\section{Stable bundles on $S$}
\label{stablebundles}

In this section we note some facts about stable bundles on the surface\linebreak
\mbox{$S=\Sigma\times\cpone$}, and in particular about the set 
$\calz_{(\cald,c)}$
defined in Section \ref{defnot}.

Let $\sigma$ denote the \poincare dual of the homology class of $S$ 
represented by \mbox{$\Sigma\times\{{\rm point}\}$}, 
and $f$ the dual of the class represented by 
\mbox{$\{{\rm point}\}\times\cpone$}.  Then $H^2(S,\rr)=H^{(1,1)}(S,\rr)$ is
generated by the classes $\sigma$ and $f$.  The intersection pairing
gives
$$\sigma\cdot\sigma=f\cdot f=0,\quad \sigma\cdot f=1.$$
The \kahler cone $\calk(S)$ is the
quadrant $\{a\sigma +bf\,|\,a,b>0\}$; any such class is represented by the
associated $(1,1)$--form of some \kahler metric on $S$.

We now recall the definition of stability in the sense of Mumford and
Takemoto.  Let $X$ be a complex manifold of dimension $n$.
Let $[\omega]\in H^{(1,1)}(X,\rr)$ be an element of  
$\calk(X)$.
For any torsion-free coherent sheaf $\xi$ over $X$ we set
$$\mu(\xi)=c_1(\xi)\cdot [\omega]^{n-1}/{\rm rank}\,\xi.$$
The sheaf $\xi$ is {\it $[\omega]$--slope stable} (resp. {\it semistable})
if for any proper coherent
subsheaf $\zeta$ of $\xi$ we have 
$$\mu(\zeta)<\mu(\xi)
\qquad(\mbox{resp.}\quad\mu(\zeta)\le\mu(\xi)).$$
For rank 2 sheaves on a surface, stability with respect to two different 
elements of the \kahler cone
is equivalent if and only if they represent points in the same
chamber (see \cite[page 142]{f}, \cite{q}).  The chamber structure depends on
the Chern classes of the sheaves in question.  Any choice of $c_1, c_2$
determines a finite set of walls, in the \kahler cone $\calk$,
as follows: any $\zeta\in H^{(1,1)}(X,\rr)$ satisfying
$$
\left\{
\begin{array}{c}
\zeta=c_1\pmod{2}\vrule width 0pt depth 5pt\\
c_1^2-4c_2\le\zeta^2<0
\end{array}
\right\},
$$
determines a wall 
$\zeta^\bot=\{x\in \calk\,|\,x\cdot\zeta=0\}$.  Denote by  $\calw$ the
union of all such walls.

The chambers are then the connected components of $\calk - \calw$.
In our case
the \kahler cone is just the first quadrant in $H^{(1,1)}(S,\rr)\cong\rr^2$,
and the walls are rays with positive slope.
See Figure \ref{chambers} for the chamber structure in $\calk(S)$ for some
choices of $(c_1,c_2)$.  Note that as $4c_2-c_1^2$ increases, so does the 
number of walls.

\setlength{\unitlength}{1cm}
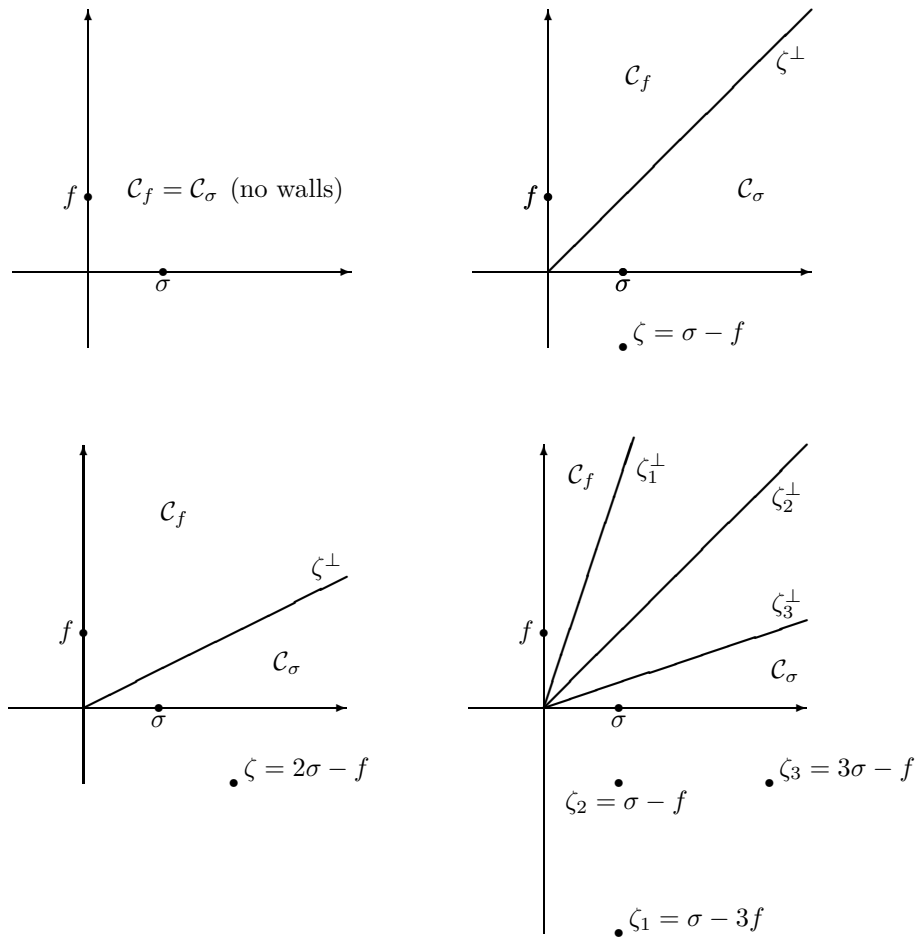
\begin{figure}[h]\small
\begin{center}
\begin{minipage}[t]{6cm}
\begin{picture}(4.5,4.5)
\put(0,1){\vector(1,0){4.5}}
\put(1,0){\vector(0,1){4.5}}
\put(2,1){\circle*{0.1}}
\put(1,2){\circle*{0.1}}
\put(2,0.9){\makebox(0,0)[t]{$\sigma$}}
\put(0.9,2){\makebox(0,0)[r]{$f$}}
\put(1.5,2){$\calc_f=\calc_\sigma$ (no walls)}
\end{picture}

\end{minipage}
\begin{minipage}[t]{6cm}
\begin{picture}(4.5,4.5)
\put(0,1){\vector(1,0){4.5}}
\put(1,0){\vector(0,1){4.5}}
\put(2,1){\circle*{0.1}}
\put(1,2){\circle*{0.1}}
\put(2,0.9){\makebox(0,0)[t]{$\sigma$}}
\put(0.9,2){\makebox(0,0)[r]{$f$}}
\put(2,1){\circle*{0.1}}
\put(1,2){\circle*{0.1}}
\put(2,0.9){\makebox(0,0)[t]{$\sigma$}}
\put(0.9,2){\makebox(0,0)[r]{$f$}}
\thicklines
\put(1,1){\line(1,1){3.5}}
\put(4,4){\makebox(0,0)[tl]{$\zeta^\bot$}}
\put(2,0){\circle*{0.1}}
\put(2.1,0.1){$\zeta=\sigma-f$}
\put(2,3.5){$\calc_f$}
\put(3.5,2){$\calc_\sigma$}

\end{picture}
\end{minipage}

\bigskip
\bigskip
\bigskip
\begin{minipage}[t]{6cm}
\begin{picture}(4.5,6.5)
\put(0,3){\vector(1,0){4.5}}
\put(1,2){\vector(0,1){4.5}}
\put(2,3){\circle*{0.1}}
\put(1,4){\circle*{0.1}}
\put(2,2.9){\makebox(0,0)[t]{$\sigma$}}
\put(0.9,4){\makebox(0,0)[r]{$f$}}
\thicklines
\put(1,3){\line(2,1){3.5}}
\put(4,4.7){\makebox(0,0)[bl]{$\zeta^\bot$}}
\put(3,2){\circle*{0.1}}
\put(3.1,2.1){$\zeta=2\sigma-f$}
\put(2,5.5){$\calc_f$}
\put(3.5,3.5){$\calc_\sigma$}
\end{picture}
\end{minipage}
\begin{minipage}[t]{6cm}
\begin{picture}(4.5,6.5)
\put(0,3){\vector(1,0){4.5}}
\put(1,0){\vector(0,1){6.5}}
\put(2,3){\circle*{0.1}}
\put(1,4){\circle*{0.1}}
\put(2,2.9){\makebox(0,0)[t]{$\sigma$}}
\put(0.9,4){\makebox(0,0)[r]{$f$}}
\thicklines
\put(1,3){\line(1,3){1.2}}
\put(2.2,6){\makebox(0,0)[bl]{$\zeta_1^\bot$}}
\put(1,3){\line(1,1){3.5}}
\put(4,6){\makebox(0,0)[tl]{$\zeta_2^\bot$}}
\put(1,3){\line(3,1){3.5}}
\put(4,4.2){\makebox(0,0)[bl]{$\zeta_3^\bot$}}
\put(2,0){\circle*{0.1}}
\put(2.1,0.1){$\zeta_1=\sigma-3f$}
\put(2,2){\circle*{0.1}}
\put(2.1,1.9){\makebox(0,0)[t]{$\zeta_2=\sigma-f$}}
\put(4,2){\circle*{0.1}}
\put(4.1,2.1){$\zeta_3=3\sigma-f$}
\put(1.3,6){$\calc_f$}
\put(4,3.4){$\calc_\sigma$}

\end{picture}
\end{minipage}
\bigskip

\caption{\label{chambers}
Chamber Structure for 
$(c_1,c_2)=(f,0),\,(\sigma+f,1),\,(f,1),\,(\sigma+f,2)$}
\end{center}
\end{figure}

We will be concerned with the chamber $\calc_\sigma$, whose 
boundary contains the point $\sigma$.  A bundle or sheaf is said to be
$\calc_\sigma$--stable if it is $[\omega]$--stable for any 
$[\omega]\in\calc_\sigma$.  Observe that for any fixed choice of $(c_1,c_2)$,
the point $n\sigma+f$ will be in $\calc_\sigma$ for large enough $n$.

Fix some line bundle $\cald$ on $S$ with $\sigma\cdot c_1(\cald)=1$
and some integer $c$.  Then we denote by
$\calm_{(\cald,c)}$ the set of isomorphism classes of 
$\calc_\sigma$--stable 
bundles with determinant and second Chern class given by $\cald$ and $c$ 
respectively.

For any choice of $\cald,c$, let ${\tilde E}\to S$ be a smooth rank two complex
vector bundle with $c_1({\tilde E})=c_1(\cald)$ and $c_2({\tilde E})=c$.
Let $\cala^{(0,1)}({\tilde E})$ denote the space of all 
integrable $(0,1)$ connections 
on ${\tilde E}$ (with the $C^\infty$ topology) and let 
$\calg^c$ denote the group
of complex linear bundle automorphisms of ${\tilde E}$.  Then we topologise
$\calm_{(\cald,c)}$ as a subspace of $\cala^{(0,1)}({\tilde E})$/$\calg^c$.  
That is
to say, a sequence of elements $[\cale_n]$ of $\calm_{(\cald,c)}$ 
converge to a limit $[\cale_\infty]$
if and only if some sequence of $(0,1)$ connections $\dbar_n$ on ${\tilde E}$
representing $[\cale_n]$ converge in $C^\infty$ to a limit $\dbar_\infty$
which represents $[\cale_\infty]$.

We will often use the term ``$\dbar$--operator''
to describe integrable $(0,1)$ connections.

The following facts about $\calc_\sigma$--stable bundles follow immediately
from \cite[Chapter 6, Theorem 5]{f} and the fact that a rank 2, odd degree
semistable bundle on $\Sigma$ must be stable.

\begin{proposition}
Let $\cale\to S$ be a rank 2 holomorphic vector bundle with 
$\det(\cale)=\cald$ and $c_2(\cale)=c$.
\begin{itemize} 
\item If $\cale$ is $\calc_\sigma$--stable then its restriction to
$\Sigma\times\{z\}$ is stable for all but finitely many $z\in\cpone$.
\item If the restriction of $\cale$ to $\Sigma\times\{z\}$ is stable for 
some $z\in\cpone$ then $\cale$ is $\calc_\sigma$--stable.
\end{itemize}
\end{proposition}

Note it follows that $\calz_{(\cald,c)}$ is a subset of
$\calm_{(\cald,c)}$.

Given a bundle $\cale$ with $[\cale]\in\calm_{(\cald,c)}$ which is unstable on
$\Sigma\times\{z_i\}$ for some finite set $z_1,\dots,z_k\in\cpone$, one may
apply a finite sequence of elementary modifications to obtain a bundle $\cale'$
which is stable on $\Sigma\times\{z\}$ for all $z\in\cpone$, and such that
$$\cale|_{\Sigma\times (\cpone-\cup z_i)}
\cong\cale'|_{\Sigma\times (\cpone-\cup z_i)} .$$
(See \cite[pages 41, 98, 148]{f}.)

For any $z\in\cpone$,
$$\det(\cale|_{\Sigma\times\{z\}})=\cald|_{\Sigma\times\{z\}}.$$
This gives the same element of ${\rm Pic}\,\Sigma$ for all $z$.  Thus for each
$z\in\cpone-\cup z_i$, the restriction $\cale|_{\Sigma\times\{z\}}$
determines a point in the moduli space $\caln(\Sigma)$.
This implies (see Section \ref{finemoduli}) 
that in fact $[\cale]$ determines a holomorphic map
$$f_\cale\co \cpone-\cup z_i\to\caln(\Sigma).$$
Also $\cale'$ determines an extension of this holomorphic map to all of 
$\cpone$, which we still denote by $f_\cale$.

The following Proposition shows that in fact
$[\cale]\in\calm_{(\cald,c)}$ is completely determined by the map $f_\cale$
and by local information at the curves $\Sigma\times\{z_i\}$ on which
$\cale$ has unstable restriction.

\begin{proposition} 
\label{triples}
A $\calc_\sigma$--stable bundle $\cale$ which satisfies
$$
\sigma\cdot c_1(\cale)=1,\quad
f\cdot c_1(\cale)=0\,\,\mbox{or}\,\, 1
$$
determines, and is determined up to isomorphism by, a triple\newline 
$(f_\cale,\cup_{i=1}^k z_i,\cup_{i=1}^k [\cale_i])$,
where 
\begin{itemize}
\item $f_\cale$ is a holomorphic map from $\cpone$ to $\caln(\Sigma)$,
\item $z_1,\dots,z_k$ are distinct (unordered) points in $\cpone$ $(k\ge0)$, 
and 
\item $[\cale_i]$ is an $f_\cale$--compatible holomorphic vector bundle 
over a germ of a neighbourhood of $\Sigma\times \{z_i\}$,
ie, an equivalence class of bundles $\cale_i\to\Sigma\times D_i$ where $D_i$ 
is a small disk neighbourhood of $z_i$ in $\cpone$, 
$\cale_i|_{\Sigma\times\{z_i\}}$ is unstable, and
$$\cale_i|_{\Sigma\times (D_i-z_i)}\cong 
({\rm Id}\times f_\cale|_{D_i-z_i})^*\calu.$$
Two such bundles $\cale_i\to\Sigma\times D_i$ and 
$\cale_i'\to\Sigma\times D_i'$ are equivalent if they are isomorphic on
$\Sigma\times (D_i\cap D_i')$.
\end{itemize}
\end{proposition}
\proof
It is clear from the previous discussion that a $\calc_\sigma$--stable bundle 
determines such a triple.  Conversely, suppose we are given
$(f_\cale,\cup_{i=1}^k z_i,\cup_{i=1}^k [\cale_i])$.
We construct a bundle $\cale$ over $S$ by gluing together
the bundle $({\rm Id}\times f_\cale|_{\cpone-\cup z_i})^*\calu$ on
$\Sigma\times(\cpone-\cup z_i)$ and the bundles $\cale_i$ on 
$\Sigma\times D_i$.  The gluing isomorphisms are well-defined up to
multiplication by a $\cc^*$--valued function on each $D_i-z_i$.  Thus the
isomorphism class of the glued-up bundle $\cale$ is well-defined up to 
tensoring with a line bundle pulled back from $\cpone$. The requirement that
$$f\cdot\det(\cale)=0\,\,\mbox{or}\,\, 1$$
determines $[\cale]$ uniquely.
\endpf

Note that a triple $(f_\cale,\cup_{i=1}^k z_i,\cup_{i=1}^k [\cale_i])$
corresponding to an element of $\calm_{(\cald,c)}$ is in the subset 
$\calz_{(\cald,c)}$ if and only if neither of the points $\{0,\infty\}$ is
included in $\cup_{i=1}^k z_i$.  Also note that reordering the points
$z_1,\dots,z_n$ and, correspondingly, the data
$[\cale_1],\dots,[\cale_n]$ does not change the bundle $\cale$. 

Let $\calm'_{(\cald,c)}$ denote the set of all triples for which the
associated stable bundle has second Chern class $c$ and determinant $\cald$.
We topologise $\calm'_{(\cald,c)}$ as follows: a sequence of triples
$(f_n,\cup_{i=1}^{k_n} z_{i;n},\cup_{i=1}^{k_n} [\cale_{i;n}])$
is converging to a limit 
$(f_\cale,\cup_{i=1}^k z_i,\cup_{i=1}^k [\cale_i])$ if and
only if
\begin{itemize}
\item $f_n \to f_\cale$ in $C^\infty$ on compact subsets of 
$\cpone-\cup_{i=1}^k z_i$
\item for {\it any} open neighbourhood $D$ of $\cup_{i=1}^k z_i$, 
there exists $N$ such that \linebreak
$\cup_{i=1}^{k_n} z_{i;n}\subset D\quad \forall\quad n\ge N$
\item for {\it some} open neighbourhood $D$ of $\cup_{i=1}^k z_i$, the 
isomorphism classes of bundles
determined by $f_n$ and $\cup_{i=1}^{k_n} [\cale_{i;n}]$ converge
on $\Sigma\times D$ to that determined by $f_\cale$ and
$\cup_{i=1}^k [\cale_i]$.  (Convergence here means that a choice of
representatives of the isomorphism classes converges in $C^\infty$ on 
compact sets to a representative of the limit.)
\end{itemize}

Proposition \ref{triples} gives a correspondence between $\calm_{(\cald,c)}$
and $\calm'_{(\cald,c)}$.
We will see in Section \ref{convergence} that in fact this is a
homeomorphism.

\section{From instantons to stable bundles}
\label{instantontostable}
Let $e\in4\pi^2\zz$ and let $\cald,c$ be determined 
by $e$ as in (\ref{det}),(\ref{c2}).
In this section we define a function
$$\Psi\co \calm_e/\cali\to\calz_{(\cald,c)}.$$
Let $A$ be a $U(2)$--instanton with energy $e$.  
Denote by $\dbar_A$ the associated
$(0,1)$ operator on \mbox{$E$--valued} forms.
We will show that the holomorphic bundle $(E,\dbar_A)$
has a unique extension $\cale$ to the compact surface $S$, 
up to isomorphism, which satisfies
\begin{itemize}
\item $\det \cale = \cald$
\item $c_2(\cale) = c$
\item $\cale|_{\Sigma_0},\cale|_{\Sigma_\infty}$ are stable.
\end{itemize}

If $g$ is a gauge transformation then 
$(E,\dbar_A)\cong(E,\dbar_{g^*A})$ as holomorphic bundles.  This shows that
the extension gives a function $\Psi$ on $\calm_e$, defined by
$$\Psi([A])=[\cale].$$
Proposition \ref{involution} shows that $\Psi$ descends to $\calm_e/\cali$.

We prove two propositions in this section.
The first shows how to extend the holomorphic vector bundle $(E,\dbar_A)$
in a unique way over the compact surface $S$.
This could be established following the method
of Guo \cite{guo}; however the product structure of $Y$ and the fact that
$\caln(\Sigma)$ is a fine moduli space give a much easier proof in our case.
The second proposition confirms that the determinant line bundle and the
second Chern class of the extended bundle are given by $\cald$ and $c$, 
respectively.

\begin{proposition}
\label{extension}
Let $A$ be a (finite energy) $U(2)$--instanton on 
$E\to \Sigma \times S^1 \times \rr$.  Then the bundle 
$E$ with the holomorphic structure determined by $A$
extends uniquely to a holomorphic bundle $\cale$ over $S=\Sigma \times \cpone$
whose restriction to $\Sigma_0$ and $\Sigma_\infty$ is stable and whose
restriction to every fiber $\{{\rm point}\}\times\cpone$ has degree $0$ or $1$.
\end{proposition}

\proof 
For convenience, we will work
with the connection ${\tilde A}$ obtained from $A$ as in Lemma \ref{gaugereps}.

The proof that the bundle extends over the divisors at infinity is the
same for each end.
We consider the end $\Sigma_0$ (corresponding to $t=-\infty$).

Fixing a bundle isomorphism between $E$ and $E_Y \times \rr$ as in
Section \ref{instantons}, ${\tilde A}$ gives an
associated path $[{\tilde A}_t]$ of gauge-equivalence classes of connections on
$E_Y$.  This path converges to projectively flat limits 
$[{\tilde A}_{\pm\infty}]$ as $t\to\pm\infty$.
By applying a bundle automorphism of $E_Y$, we may arrange that 
${\tilde A}_{-\infty}$ is pulled back from a projectively flat connection
$a_{-\infty}$ on $\Sigma$.
Since $[{\tilde A}_t] \to [{\tilde A}_{-\infty}]$ as $t \to -\infty$, 
it follows that the restriction
$[{\tilde A}_{(\theta,t)}]$ to $E|_{\Sigma \times (\theta,t)}$ is converging to
$[a_{-\infty}]$.  Thus also the restriction of the 
holomorphic structure determined
by $[{\tilde A}]$ to $E|_{\Sigma \times (\theta,t)}$ is converging to that 
given by $[a_{-\infty}]$.

Now $a_{-\infty}$ is a projectively flat connection on $E_\Sigma$; 
therefore the
holomorphic structure it determines is stable.  Stability is an open
condition on the space of equivalence classes of holomorphic structures.
Thus for some $T<<0$, we must have that $[{\tilde A}_{(\theta,t)}]$ 
determines a stable holomorphic structure for each $t<T$.  Or in other words
$(E,\dbar_{\tilde A})$ gives a holomorphic bundle on the product
of $\Sigma$ times the punctured disk 
$D_0{}\!^*=\{(\theta,t):t<T\}=\{z:0<|z|<R=e^{T}\}$, whose restriction
to each ${\Sigma \times \{z\}}$ is stable with determinant $\call$, i.e an 
element of the moduli space $\caln(\Sigma)$.

Recall from Section \ref{finemoduli} that $\caln(\Sigma)$ is a 
fine moduli space.
Thus there exists a holomorphic map
$$f_0\co D_0{}\!^*\to \caln(\Sigma)$$
such that
$$(E,\dbar_{\tilde A})|_{\Sigma \times D_0{}\!^*}
\cong({\rm Id}\times f_0)^*\calu.$$
(Noting that there are no nontrivial holomorphic line bundles on $D_0{}\!^*$.)
Denote the isomorphism between them by $g_0$.  This is a section of
$GL(E)$ over $\Sigma \times D_0{}\!^*$.

The fact that $f_0(z)$ is approaching the point in $\caln(\Sigma)$ determined
by $a_{-\infty}$ as $z \to 0$ implies by removable singularities 
that $f_0$ can be extended holomorphically
to a map of the disk $D_0$ by setting $f_0(0)=[\dbar_{a_{-\infty}}]$.
Then $\cale_0=({\rm Id}\times f_0)^*\calu$ is a holomorphic bundle on 
$\Sigma \times D_0$ whose restriction to 
$\Sigma \times D_0{}\!^*$ is isomorphic to $(E,\dbar_{\tilde A})$.

In exactly the same way we can find a disk neighbourhood $D_\infty$ of
$\infty\in\cpone$ and a bundle $\cale_\infty$ on $\Sigma \times D_\infty$
whose restriction to $\Sigma \times D_\infty{}\!^*$ is isomorphic to 
$(E,\dbar_{\tilde A})$ by some $g_\infty\in\Gamma(\Sigma\times D_\infty{}\!^*,GL(E))$.

Now we can patch together these bundles using the complex gauge
transformations $g_0,g_\infty$ to form a bundle $\cale$ on
$\Sigma\times\cpone$ which extends the holomorphic bundle
$(E,\dbar_{\tilde A})$ on the cylinder.  This extension is not unique;
the universal property of $\calu\to\Sigma\times\caln(\Sigma)$ only
defines the bundles $\cale_0,\cale_\infty$ up to tensoring with a line
bundle pulled back from $D_0, D_\infty$ respectively.  Thus the
extension $\cale$ is unique up to tensoring with a line bundle pulled
back from $\cpone$. The degree of the restriction of $\cale$ to a
fiber $\{{\rm point}\}\times\cpone$ is given by the intersection
pairing $c_1(\cale)\cdot f$.  Tensoring with $\pi_2^*\calo_\cpone(a)$
will change this intersection by $2a$; thus there will be a unique
extension $\cale$ with ${\rm deg}\,\cale|_{\{{\rm
point}\}\times\cpone}=0$ or $1$.  \endpf

\begin{remark}
\label{extensionremark}
Proposition \ref{extension} does not give an extension of the
connection ${\tilde A}$ to $S$, or even an extension of the unitary
structure on $E$.  We just see that the holomorphic bundle extends
over $S$.  We will see from the proof of the next Proposition that in
fact the metric extends {\em continuously} to the extended holomorphic
bundle.  This is reminiscent of Guo's results \cite[I, Theorem
6.1]{guo}.

\end{remark}

\begin{proposition}
\label{det&c2}
Let $A$ be a $U(2)$--instanton on 
$E\to \Sigma \times S^1 \times \rr$ with energy $e(A)=e\in4\pi^2\zz$.  
Let $\cale$ be the extension of the holomorphic bundle $(E,\dbar_A)$ 
over $S$ given in Proposition \ref{extension}, and let $\cald,c$ be given
by {\rm (\ref{det}), (\ref{c2})} in Section \ref{defnot}.  Then the determinant
line bundle of $\cale$ is isomorphic to $\cald$ and $c_2(\cale)=c$.
\end{proposition}

\proof
We first consider the determinant line bundle $\det\cale$.  
Note that the central part
of the connection $A$ and the central parts of the connections $A_{\pm\infty}$ 
and $a_{\pm\infty}$ described in the proof of Proposition \ref{extension} are
all given by the same fixed central connection pulled back from $\Sigma$.
This connection determines a degree 1 holomorphic line bundle $\call$ on 
$\Sigma$; it follows that the restriction of $\cale$ to each fiber 
$\Sigma\times \{{\rm point}\}$ has $\call$ as its determinant line bundle.
This then implies that $\det\cale\cong\pi_1^*\call+\pi_2^*\calo_\cpone(a)$ for
some $a\in\zz$.  The
specific choice of the extension $\cale$ implies that $a=0\,\,\mbox{or}\,\,1$. 
It follows that $c_1(\cale)=a \sigma + f$, so that $c_1(\cale)^2=2a$.  
The result will now follow if we can establish the following identity:
\begin{equation}
\label{chernweil}
c_1(\cale)^2-4c_2(\cale)=-\frac{1}{2\pi^2}e(A).
\end{equation}
To see this note first that, by (\ref{chernweil})
$$c_1(\cale)^2=2a=-\frac{1}{2\pi^2}e(A)\bmod4$$
which implies that
$$
a=
\left\{\begin{array}{ll}
0&\mbox{if}\quad e(A)\in8\pi^2\zz,\\
1&\mbox{otherwise}.
\end{array}\right.
$$
This tells us that $\det\cale=\cald$ as required; it then follows immediately
from (\ref{chernweil}) that $c_2(\cale)=c$.

We will establish Equation (\ref{chernweil}) using Chern--Weil theory.
For any smooth connection $A'$ on
a rank 2 complex vector bundle $\cale'$ over a compact Riemannian 4--manifold
$X$ we have the following identity:
$$p_1({\rm ad}\,\cale')=c_1(\cale')^2-4c_2(\cale')=
-\frac{1}{2\pi^2}\int_X\tr F_{A'}^0\wedge F_{A'}^0.$$
This number is always an integer, and it only depends on the toplogical type
of the bundle $\cale'$.  Thus to prove (\ref{chernweil}) it suffices to find
a complex
bundle $\cale'$ which is topologically equivalent to $\cale$, and a 
connection $A'$
on $\cale'$ such that the integral $\int_X\tr F_{A'}^0\wedge F_{A'}^0$ 
differs from $e(A)$ by less than $2\pi^2$.  We will do this by gluing together
the connection ${\tilde A}$ on the cylinder and flat connections on the 
divisors at infinity, where ${\tilde A}$ is obtained from $A$ as in 
Lemma \ref{gaugereps}.

From the proof of Proposition \ref{extension} the bundle $\cale$ is formed by
patching together bundles on the cylinder and on neighbourhoods of the divisors
at infinity as follows:
$$\cale=\cale_0\cup_{g_0}(E,\dbar_{\tilde A})\cup_{g_\infty}\cale_\infty,$$
where $g_0, g_\infty$ are sections of $GL(E)$ over 
$\Sigma\times D_0{}\!^*, \Sigma\times D_\infty{}\!^*$ respectively.

Denote by $E_0, E_\infty$ the smooth hermitian vector bundles underlying
$\cale_0,\cale_\infty$.  Both of these vector bundles are smoothly
bundle isomorphic to the pullback of the degree 1 hermitian vector bundle
$E_\Sigma$ on $\Sigma$.

We again restrict our attention to one end.  Without loss of generality, we 
consider the end corresponding to $z=0$ in $\cpone$ (ie $t\to-\infty$).  
After a unitary gauge change, we may suppose that ${\tilde A}$ is
converging exponentially fast to a limit $\pi_\Sigma^*a_{-\infty}$,
where $a_{-\infty}$ is a projectively flat connection on $\Sigma$.
(See Lemma \ref{gaugereps}.)

Let $\dbar_{\tilde A}$ denote the $(0,1)$ operator associated to 
${\tilde A}$, and let $\dbar_{a_{-\infty}}$ be
the $(0,1)$ operator on $\Sigma$ associated to $a_{-\infty}$.
It follows that
the restriction $\dbar_{\tilde A}|_{\Sigma_z}$ is converging to 
$\dbar_{a_{-\infty}}$
 as $z\to 0$. (Here $\Sigma_z$ denotes the curve 
$\Sigma\times\{z\}\subset\Sigma\times\cpone$.)

Denote by $\dbar_B$ the $(0,1)$ operator on $E_0\to\Sigma\times D_0$ which
gives the holomorphic structure on $\cale_0$.  Then
$$\dbar_B=g_0^{-1}\circ\dbar_{\tilde A} \circ g_0.$$
The restriction of $\dbar_B$ to $\Sigma_0$ is by construction isomorphic to
$\dbar_{a_{-\infty}}$.  
Thus after possibly changing $\dbar_B$ and $g_0$ by a complex
gauge automorphism pulled back from $\Sigma$ we may assume that
$\dbar_B|_{\Sigma_0}=\dbar_{a_{-\infty}}$.  Thus we have
$$\dbar_{\tilde A}|_{\Sigma_z}\to\dbar_{a_{-\infty}}\quad\mbox{as}\ z\to 0,$$
and also
$$g_0^{-1}\circ\dbar_{\tilde A}\circ g_0|_{\Sigma_z}
\to\dbar_{a_{-\infty}}\quad\mbox{as}\ z\to 0.$$
We may also assume that the central part of $\dbar_B$ is pulled back from
$\Sigma$ and in fact that it agrees on $\Sigma\times D_0{}\!^*$ 
with the central part of $\dbar_{\tilde A}$.  It then follows that
$\det(g_0)$ is a holomorphic function from $D_0{}\!^*$ to $\cc^*$.

\begin{claim}
\label{even}
The winding number of $\det(g_0)$ is even.
\end{claim}

Assume for now that the claim is true.

Then we may replace $g_0$ by $h=\frac{1}{\sqrt{\det{g_0}}}g_0$; this is a
section of $SL(E)$ over $\Sigma\times D_0{}\!^*$, that is a
complex gauge automorphism whose determinant is the constant function 1.
Since $\frac{1}{\sqrt{\det{g_0}}}$ is a holomorphic function, $h$ also
satisfies the equation
$$\dbar_B=h^{-1}\circ\dbar_{\tilde A}\circ h.$$
Also as $z\to 0$, $h|_{\Sigma_z}$ must be converging to an element in the
stabiliser of $\dbar_{a_{-\infty}}$.  That is to say,
$$h|_{\Sigma_z}\to \pm{\bf 1}\quad\mbox{as}\ z\to 0.$$
(The stabiliser consists of constants times the identity, since stable
bundles are simple \cite[page 88]{f}.)
We may suppose it is converging to $+{\bf 1}$ (multiply $h$ by $-1$ if 
necessary).

It now follows that the hermitian metric extends, at
least continuously, to the extension over $\Sigma_0$ given by $h$.\\

To establish Equation (\ref{chernweil}) we will modify the extended bundle.  
The fact that $h$ has a 
continuous extension to $\Sigma\times D_0$ with $h|_{\Sigma_0}={\bf 1}$ 
in fact shows that $h$ is continuously homotopic to the constant identity gauge
automorphism.  Indeed a homotopy is given by
$$H_s(z)=h(z-s z).$$
(Using the identification
$$\Gamma(\Sigma\times D_0{}\!^*,SL(E))=
{\rm Maps}\,(D_0{}\!^*,\Gamma(\Sigma,SL(E_\Sigma)))\,.)$$
Thus we may replace $h$ by the identity, at each end, to get
a topologically equivalent smooth bundle $\cale'$ over $\Sigma\times\cpone$.
Now the metric clearly extends to a smooth hermitian metric on this bundle,
and we can produce a smooth unitary connection
$A'$ on $\cale'$ by ``damping down''
 the original connection $A$ to the limiting connections
$\pi_\Sigma^*a_{\pm\infty}$ on each end of the cylinder.  This is a standard
technique.

Start by choosing some $T>>0$.  
It follows from Theorem \ref{thmB} and our choice of gauge representatives for
$[{\tilde A}]$ on the ends of the cylinder that we can write
$${\tilde A}=\pi_\Sigma^*a_{\pm\infty}+\alpha_\pm,$$
with 
\begin{equation}
\label{alphadecay}
||\alpha_\pm||_{L^2_2(Y\times [\pm T-\frac12,\pm T+\frac12])}
\le C e^{-\delta|T|}.
\end{equation}
We choose a smooth function $\beta_+\co \rr\to\rr$ with
$$\beta_+(t) = \left\{ \begin{array}{ll}
1 & \mbox{if} \quad t > T + \frac12\vrule width 0pt depth 10pt\\
0 & \mbox{if} \quad t < T - \frac12\,.\\
\end{array} \right.$$
Also define $\beta_-(t)=\beta_+(-t)$.
Now the connection
$$A'={\tilde A}-\beta_\pm\cdot\alpha_\pm$$
extends to the bundle $\cale'$, since it is equal to 
$\pi_\Sigma^*a_{\pm\infty}$ on the ends.

Computing the curvature, we get
\begin{eqnarray*}
F_{A'}=F_{\tilde A}
&\!\!\!-&\!\!\!
\beta_+\cdot d_{\tilde A}\alpha_+-d\beta_+\wedge\alpha_+
+\beta_+{}\!^2\cdot\alpha_+\wedge\alpha_+\\
&\!\!\!-&\!\!\!
\beta_-\cdot d_{\tilde A}\alpha_--d\beta_-\wedge\alpha_-
+\beta_-{}\!^2\cdot\alpha_-\wedge\alpha_-.
\end{eqnarray*}

By making $T$ large enough, it follows from (\ref{alphadecay}) and the usual
Sobolev multiplication theorems that the integrals
$\int_S\tr F_{A'}^0\wedge F_{A'}^0$ and 
$\int_{\Sigma\times S^1\times \rr}\tr F_{\tilde A}^0\wedge F_{\tilde A}^0=e(A)$
may be made arbitrarily close.  Since both are integer multiples of $4\pi^2$,
they must in fact be the same, as required.\endpf

\proof[Proof of Claim \ref{even}]  The connection
$\cali({\tilde A})$ (see Lemma \ref{involution}) is the image of ${\tilde A}$
under a complex gauge transformation whose determinant has winding number 1 
around $z=0$.
Suppose that the winding number of $\det(g_0)$ is odd.  Then we may write
$$\dbar_B=(g_0')^{-1}\circ\dbar_{\cali({\tilde A})} \circ g_0',$$
with $\det(g_0')$ even.  Then as above we may divide by the square root of
the determinant to obtain a section $h$ of $SL(E)$ over 
$\Sigma\times D_0{}\!^*$ satisfying
$$\dbar_B=(h)^{-1}\circ\dbar_{\cali({\tilde A})} \circ h.$$
We will show that in fact this is not possible by comparing the asymptotic
limits of the operators
$\dbar_B$ and $\dbar_{\cali({\tilde A})}$ on the end of the cylinder.

Since $\dbar_B$ extends smoothly over $\Sigma\times D_0$, 
and its restriction to 
$\Sigma_0$ is $\dbar_{a_{-\infty}}$ we have
$$\lim_{t\to-\infty}\dbar_B|_{\Sigma_{(\theta,t)}}=\dbar_{a_{-\infty}},$$
with convergence in $C^\infty$ and hence also in $L^2_1$, uniformly in 
$\theta$.

Theorem \ref{thmB} and the trace theorem for Sobolev spaces (see for example
\cite{bw}) imply that
$$\lim_{t\to-\infty}\dbar_{\cali({\tilde A})}|_{\Sigma_{(\theta,t)}}=
\dbar_{\cali(\pi_\Sigma^*a_{-\infty})}|_{\Sigma_{(\theta,t)}},$$
where the limit is taken in the $L^2_1$ topology.  By the continuity property 
of the trace operator, this convergence is also uniform in $\theta$.
Note that the restriction 
$\dbar_{\cali(\pi_\Sigma^*a_{-\infty})}|_{\Sigma_{(\theta,t)}}$
is independent of $t$; in fact it follows from the definition of $\cali$
that
$$\dbar_{\cali(\pi_\Sigma^*a_{-\infty})}|_{\Sigma_{(\theta,t)}}=
\left(\begin{array}{cc}
e^\frac{i\theta}2 & 0\\
0 & e^{-\frac{i\theta}2} \end{array} \right)^*\dbar_{a}.$$
The above uniform convergences, together with the slice condition (see Section
\ref{slices}) imply that
$$\lim_{t\to-\infty}h|_{\Sigma_{(\theta,t)}}=
\pm\left(\begin{array}{cc}
e^\frac{i\theta}2 & 0\\
0 & e^{-\frac{i\theta}2} \end{array} \right),$$
with convergence in $L^2_2$ and hence in $C^0$, uniformly in $\theta$.
This is a contradiction: it
is not possible to have a smooth (or even continuous) complex gauge
transformation $h$ on the cylinder $Y\times S^1$ converging uniformly to
a discontinuous limit on $Y$.
\endpf

\section{Injectivity of the function $\Psi$}
\label{injectivity}

Proposition \ref{extension} defines a function
$$\Psi\co \calm_e/\cali\to\calz_{(\cald,c)}.$$
In this section we establish that the map $\Psi$ is injective.

\begin{proposition}
\label{injective}
$\Psi$ is injective.  That is, if $[A],[A']\in\calm_e$ determine isomorphic
holomorphic structures and both $[A_t]$ and $[A_t']$ converge to limits in
$\rplus$ as $t\to\infty$, then $[A]=[A']$. 
\end{proposition}

(Recall that the involution $\cali$ switches the components $\rplus,\rminus$
of projectively flat connections over $Y$; fixing $[A_\infty]\in\rplus$
chooses between $[A]$ and $\cali([A])$.)

The proof of Proposition \ref{injective} uses the maximum principle applied to 
a function measuring
the distance between two metrics, following Donaldson \cite{d1,d3}.  
(Much of the proof is taken directly from \cite{d3}).

If $H$ and $K$ are Hermitian metrics on $E$ 
(which may be regarded as isomorphisms
from $E$ to the dual bundle $E^*$) 
then $\eta=H^{-1}K$ is a section of the bundle of 
endomorphisms ${\rm End}\,E$ which is self-adjoint
with respect to either metric.
\begin{definition}
For any two hermitian metrics $H,K$ on $E$
set
$$\sigma(H,K)={\rm Tr}(H^{-1}K)+{\rm Tr}(K^{-1}H)-4 
\quad\in C^\infty(\Sigma\times S^1\times\rr).$$
\end{definition}

This has the property that $\sigma(H,K)\ge0$ with equality if and only if
$H=K$.

Let $H=\langle\,\cdot\,,\cdot\,\rangle$ represent the metric on $E$ 
which is pulled back from that on 
$E_\Sigma$, and with respect to which $A, A'$ are unitary.  By hypothesis
$$\dbar_{A'}=g^{-1}\circ\dbar_A\circ g$$
for some complex gauge transformation $g$, whose determinant we may require
to be identically 1.
Let $K$ be the Hermitian metric on $E$ given by
$$K=\langle\, g\,\cdot\,,g\,\cdot\,\rangle=
\langle\,\cdot\,,g^*g\,\cdot\,\rangle=
\langle\,\cdot\,,\eta\,\cdot\,\rangle.$$
We will show that in fact $g$ is a unitary gauge transformation by
demonstrating that $H=K$.

Let $A_H$ be the unique connection determined by the metric $H$ and the
holomorphic structure $\dbar_A$, and let $A_K$ be the connection
determined by $K$ and $\dbar_A$.  Then in fact $A_H=A$, and 
$A_K=g\circ A'\circ g^{-1}$.  Thus $A_H$ is a
$U(2)$--instanton, while $A_K$ is a finite-energy, projectively ASD connection
with the same central part as $A_H$, but it is unitary with respect to the
metric $K$ rather than $H$.

The curvature of these connections is related by
\begin{equation}
\label{HKcurvature}
F_K=F_H+\dbar_H(\eta^{-1}\partial_H\eta).
\end{equation}
For an integrable connection $A$ in a vector bundle over a compact
\kahler surface, it is not hard to show that the condition that $A$ is a
$U(2)$--instanton (ie\ ${\rm Tr}F_A$ harmonic, $F_A^0$ ASD) is equivalent
to the Hermitian--Einstein condition 
$$\hat{F_A}=\lambda\cdot{\bf 1}$$
for constant scalar $\lambda$.  (Here $\hat{F_A}=\Lambda F_A$, where
$\Lambda$ is the adjoint of wedging with the \kahler form, and ${\bf 1}$
is the identity endomorphism of the vector bundle in question).

Generalising to a possibly noncompact \kahler surface, one finds the
equivalence
$$\left\{ \begin{array}{c}F_A^0\quad\mbox{ASD}\\
{\rm Tr}F_A\quad\mbox{harmonic}\end{array}\right\}
\Leftrightarrow\hat{F_A}=f\cdot{\bf 1},$$
where $f$ is a harmonic function.  In our definition of
$U(2)$--instantons on $\Sigma\times S^1\times\rr$, we require that 
${\rm Tr}F_A$ is both harmonic and pulled back from $\Sigma$.  This
means that $A_H$ and $A_K$ both satisfy the Hermitian--Einstein condition
$\hat{F_A}=\lambda\cdot{\bf  1}$, with $\lambda$ a fixed constant
depending on the first Chern class of the restriction of the bundle to 
$\Sigma\times \{{\rm point}\}$.

\begin{lemma}
\label{subharmonic}
$\Delta\sigma(H,K)\le0.$
\end{lemma}

\proof 
Applying $i\Lambda$ to $(\ref{HKcurvature})$ and taking the trace yields
\begin{eqnarray}
0 &=& {\rm Tr}\,i\Lambda\dbar_H(\eta^{-1}\partial_H\eta)\nonumber\\
&=& {\rm Tr}\,i\Lambda(\eta^{-1}\dbar_H\partial_H\eta+
\eta^{-1}\dbar_H\eta\,\eta^{-1}\wedge\partial_H\eta).
\label{eqn2}
\end{eqnarray}

Now we use the Weitzenb\"{o}ck formula (see \cite[page 212]{dk}): for any
connection $A$ on $E$,
$$\partial_A^*\partial_A=\frac12\nabla_A^*\nabla_A-i{\hat F_A}$$
on $\Omega^0({\rm End}\,E)$, with $i{\hat F_A}$ acting by the adjoint
action.  Here $\nabla_A=\partial_A+\dbar_A$ is the covariant derivative on
${\rm End}\,E$ associated to $A$.
We also have the following first order \kahler identity on
$\Omega^{1,0}({\rm End}\,E)$:
$$\partial_A^*=-i\Lambda\dbar_A.$$
Applying these identities as well as the Hermitian--Einstein condition to
$(\ref{eqn2})$ yields
$$\Delta{\rm Tr}\,\eta={\rm Tr}\,\nabla_H^*\nabla_H\eta=
2i\Lambda{\rm Tr}(\dbar_H\eta\,\eta^{-1}\wedge\partial_H\eta).$$
Choose a frame for $E$ over any given point $p$ which is unitary 
with respect to $H$
and in which $\eta$ is diagonal with eigenvalues $\lambda_a$, and let 
$\pi_{ab}$ be the matrix entries of $\dbar_H\eta$ at $p$.  Then 
$\partial_H\eta$ has entries ${\bar \pi}_{ba}$ and
$$i\Lambda{\rm Tr}(\dbar_H\eta\,\eta^{-1}\wedge\partial_H\eta)
=i\Lambda\sum_{a,b}\lambda_a^{-1}\pi_{ab}\wedge{\bar \pi}_{ab}
=-\sum_{a,b}\lambda_a^{-1}|\pi_{ab}|^2\le0.$$
(Using the fact that for a $(1,0)$ form $\phi$, 
$i\Lambda(\phi\wedge{\bar \phi})=-|\phi|^2$.)

Thus $\Delta{\rm Tr}\,(H^{-1}K)\le0$.  Interchanging $H,K$ we see that
${\rm Tr}\,(K^{-1}H)$ is likewise sub-harmonic, and hence also
$\sigma(H,K)$.\endpf

The maximum principle (see for example \cite{gt}) states that a subharmonic
function on a compact domain attains its maximum value on the boundary.  In
order to apply this we need to see what happens to $\sigma(H,K)$ on the
ends of the cylinder $\Sigma\times S^1\times\rr$.

\begin{lemma}
\label{samelimits}
$\sigma(H,K)\to0\quad\mbox{as}\quad t\to\pm\infty.$
\end{lemma}

\proof
Consider again the $U(2)$--instantons $A,A'$ which are related by the 
complex gauge transformation $g$ with determinant 1.  By Theorem \ref{thmA},
both $[A_t],[A'_t]$ converge to flat limits as $t\to\pm\infty$.  Using
Theorem \ref{thmB} and the Sobolev trace theorem (\cite{bw}), 
$A|_{\Sigma_{(\theta,t)}}$ and $A'|_{\Sigma_{(\theta,t)}}$ converge 
in $L^2_1$ to projectively flat limits on $\Sigma$ as 
$t\to\pm\infty$, uniformly in 
$\theta$.  These flat limits determine isomorphic holomorphic structures, so
by the Narasimhan--Seshadri theorem they are in the same unitary gauge orbit.
Arguing as in the proof of Claim \ref{even} (see also proof of Lemma 
\ref{uniformC0}) we find that 
$g|_{\Sigma_{(\theta,t)}}$ is converging in $C^0$ to a unitary limit for each
$\theta$ as $t\to\pm\infty$, uniformly in $\theta$.  Thus $\eta=g^*g$ and
$\eta^{-1}$ are converging to the identity in $C^0$ as $t\to\pm\infty$, from
which it follows that $\sigma(H,K)\to0$.\endpf

{\bf \noindent Proof of Proposition \ref{injective}.}
Let $H,K$ be the metrics obtained from $A,A'$ as above.
For any $\epsilon>0$, it follows from
Lemma \ref{samelimits} that we may choose $T>0$ such that 
$\sigma(H,K)<\epsilon$ for $|t|\ge T$.  Lemma \ref{subharmonic} tells us that
$\Delta\sigma(H,K)\le 0$.  It follows from the maximum principle that
$\sigma(H,K)<\epsilon$ on $\Sigma\times S^1\times [-T,T]$ and thus on all of
$\Sigma\times S^1\times\rr$.  Thus $\sigma(H,K)=0$ at each point, and $H=K$.
It follows that the gauge transformation $g$ is unitary, so that
$$[A]=[A']$$
in $\calm_e$, as required.
\endpf

\section{Some convergence results}
\label{convergence}

In this section we prove some convergence results for complex gauge
transformations.  These will be used to prove that $\Psi$ is surjective
and that $\Psi$ and its inverse are continuous.

Throughout this section $E_\Sigma\to \Sigma$ is a smooth rank two hermitian
vector bundle with degree one, and
for any subset $W\subset\cpone$, we denote by 
$$E\to \Sigma\times W$$
the pullback of $E_\Sigma$ by the projection onto the first factor.  The term
$\dbar$--operator refers to smooth $(0,1)$ connections.

In the following two lemmas, take $W$ to be an annulus in $\cpone$.
Any time we refer to an annulus contained in another
annulus, it is to be understood that the larger annulus retracts
onto the smaller one.

\begin{lemma}
\label{uniformC0}
Suppose $\{\dbar_{A_n}\}_{n=1}^\infty,\{\dbar_{B_n}\}_{n=1}^\infty,\dbar_A$ 
and $\dbar_B$ are $\dbar$--operators on \linebreak $E \to \Sigma\times W$ 
whose restrictions to $\Sigma\times\{z\}$ for each
$z\in W$ are stable with fixed determinant $\call$, and that
\begin{eqnarray*}
\dbar_{A_n}&\to&\dbar_A,\\
\dbar_{B_n}&\to&\dbar_B
\end{eqnarray*}
in $C^\infty$ on compact subsets of $\Sigma\times W$.
Suppose also that
\begin{eqnarray*}
\dbar_{B_n}&=&\phi_n^*\dbar_{A_n},\\
\dbar_B&=&\phi^*\dbar_A
\end{eqnarray*}
where $\{\phi_n\}_{n=1}^\infty,\phi$ are sections of $SL(E)$ over 
$\Sigma\times W$.
Then $\{\phi_n\}$ is uniformly bounded in $C^0(\Sigma\times W_0)$, 
where $W_0$ is a compact annulus contained in $W$.
\end{lemma}

\proof
Recall from Section \ref{slices} that there are local slices for the action
of the complex gauge group on the space of $(0,1)$ connections over $\Sigma$.

For each $z\in W$, let $\dbar_z$ denote the $(0,1)$ connection on 
$E_\Sigma\to \Sigma$ given by the restriction of $\dbar_A$ to 
$\Sigma\times \{z\}$.  Denote by
$$S_z=\dbar_z+{\rm Ker}\,\dbar_z^*$$
the slice in $\cala_0^{(0,1)}$ through $\dbar_z$.
Let $U_z$ denote the 
neighbourhood
of $\dbar_z$ in the slice with the property that the natural map
$$U_z\times_{\pm{\bf 1}}\calg_0^c\to\cala_0^{(0,1)}$$
is a diffeomorphism onto its image.
By continuity of $\dbar_A$ and by the $C^\infty$ convergence 
$\dbar_{A_n}\to\dbar_A$, it follows that for some large $N$ and some open disk
neighbourhood $D_z$ of $z$ in $W$, the restriction of $\dbar_{A_n}$ to
$\Sigma\times \{z'\}$ is in the image of $U_z\times_{\pm{\bf 1}}\calg_0$ for
all $n\ge N$ and for all $z'\in D_z$.  Then since
$$\dbar_{A_n}\to\dbar_A$$
and
$$\phi_n^*\dbar_{A_n}\to\phi^*\dbar_A,$$
it follows that for any neighbourhood $O$ of the identity in $\calg_0^c$, there
exists some $N'$ such that $(\pm\phi^{-1}\phi_n)|_{\Sigma\times \{z'\}}\in O$
for all $n\ge N'$ and for all $z'\in D_z$.

It follows that the $C^0$ norms of $\phi^{-1}\phi_n$ and hence also of
$\phi_n$ are uniformly bounded on $\Sigma\times D_z$.  Choosing a finite
cover of $\Sigma\times W_0$ by such sets $\Sigma\times D_z$ then yields the
desired result.

\begin{lemma}
Suppose $\{\dbar_{A_n}\}_{n=1}^\infty,\{\dbar_{B_n}\}_{n=1}^\infty,\dbar_A$ 
and $\dbar_B$ are $\dbar$--operators on \linebreak
\mbox{$E \to \Sigma\times W$}, and that
\begin{eqnarray*}
\dbar_{A_n}&\to&\dbar_A\\
\dbar_{B_n}&\to&\dbar_B
\end{eqnarray*}
in $C^\infty(\Sigma\times W_0)$, where $W_0$ is a compact annulus contained 
in $W$.  Suppose also that
$$\dbar_{B_n}=\phi_n^*\dbar_{A_n},$$
where $\{\phi_n\}$ are a sequence of sections of $SL(E)$ over $\Sigma\times W$
which are uniformly bounded in $C^0(\Sigma\times W_0)$.
Then there exists a subsequence of $\{\phi_n\}$ which is converging to 
some limit $\phi$ in $C^\infty(\Sigma\times W_\infty)$ with
$$\dbar_B=\phi^*\dbar_A,$$
where 
$W_\infty\subset W_0$ is a slightly smaller compact annulus.
\end{lemma}

\proof
This is a standard bootstrap argument.\endpf

The above lemmas will enable us to prove a very useful convergence result.
We will first describe the rather complicated hypotheses.

Recall from Section \ref{finemoduli} that $\calu$ is the universal bundle
over $\Sigma\times\caln(\Sigma)$.

Let $U$ and $V$ be two possibly disconnected subsets of $\cpone$ such that
$$\cpone=U\cup V,$$
and $W=U\cap V$ is a disjoint union of finitely many open annuli. 
Let $X$ be any subset
of $\Sigma\times (U-W)$.
Let $\{\dbar_{A_n}\}_{n=1}^\infty,\dbar_A$ be 
$\dbar$--operators on \mbox{$E \to \Sigma\times U$}, with the following
properties:
\begin{itemize}
\item $\dbar_{A_n}\to\dbar_A$ in $C^\infty$ on compact subsets of 
$(\Sigma\times U)-X$,
\item $(E,\dbar_{A_n})\cong({\rm Id}\times f_n)^*\calu
\quad\mbox{on}\quad\Sigma\times W$, and
\item $(E,\dbar_{A})\cong({\rm Id}\times f_\infty)^*\calu
\quad\mbox{on}\quad\Sigma\times W$,
\end{itemize}
where 
$f_n\co \cpone\to\caln(\Sigma)$
are a sequence of holomorphic maps converging in $C^\infty$ on compact
subsets of $\Sigma\times V$ to a limit 
$f_\infty$.

One may form bundles $\cale_n$ over $\Sigma\times \cpone$ by gluing
$(E,\dbar_{A_n})$ and $({\rm Id}\times f_n)^*\calu|_{\Sigma\times V}$
along the overlap region $W$, and similarly a bundle $\cale_\infty$
from $(E,\dbar_A)$ and 
\mbox{$({\rm Id}\times f_\infty)^*\calu|_{\Sigma\times V}$}.

\begin{corollary}
\label{magic}
Let $[\cale_n]$ and $[\cale_\infty]$ be as described above.  
If all of $[\cale_n]$ have isomorphic determinant line bundles and
the same second Chern class, then each $\dbar_{A_n}$ may be extended to
a $\dbar$--operator $\dbar_n$ on a smooth vector bundle 
${\tilde E}\to S$, and $\dbar_A$ may be extended to $\dbar_\infty$
on ${\tilde E}|_{S-X}$,
so that the following hold:
\begin{itemize}
\item $({\tilde E},\dbar_n)\cong\cale_n$,
\item after passing to a subsequence, $\dbar_n$ converge to a limit 
$\dbar_\infty$ in $C^\infty$ on compact 
subsets of $S-X$, and
\item $({\tilde E},\dbar_\infty)|_{S-X}\cong\cale_\infty|_{S-X}.$
\end{itemize}

In particular if $X=\emptyset$, then (after passing to a subsequence) 
$[\cale_n]\to[\cale_\infty]$ in $\calm_{(\cald,c)}$.

\end{corollary}

\proof
Denote by $E_V\to\Sigma\times V$ the pullback of $E_\Sigma\to\Sigma$.
For each $n$ we may choose $\dbar_{V_n}$ on $E_V$ so that
$$(E_V,\dbar_{V_n})\cong({\rm Id}\times f_n)^*\calu|_{\Sigma\times V}.$$
These may be chosen so that 
${\rm lim}_{n\to\infty}\dbar_{V_n}=\dbar_V$, in $C^\infty$ on compact
subsets of $\Sigma\times V$, where
$$(E_V,\dbar_V)
\cong({\rm Id}\times f_\infty)^*\calu|_{\Sigma\times V}.$$
We may also require that the central part of each $\dbar_{V_n}$ is pulled back
from $\Sigma$, and is the same for each $n$.
Identifying the bundles $E$ and $E_V$ over $\Sigma\times W$ in the obvious way,
it follows that
\begin{eqnarray*}
\dbar_{V_n}&=&\phi_n^*\dbar_{A_n},\\
\dbar_{V_\infty}&=&\phi^*\dbar_A
\end{eqnarray*}
for some sections $\phi_n,\phi$ of $GL(E)$ over $\Sigma\times W$. 
Replace the central parts of $\dbar_{A_n}$ and $\dbar_A$ by the same form
pulled back from $\Sigma$; it then follows that the determinants of
$\phi_n$ and $\phi$ will be holomorphic functions from $W$ to 
$\cc^*$.  For each $n$, the function $\det \phi_n$ will have either odd or
even winding number around each component annulus in $W$.  Pass to a
subsequence for which the parities of these winding numbers are the same for
all $n$.  

On an annulus on which the winding number is even, $\phi_n$ may
be replaced by $\phi_n'=\frac{1}{\sqrt{\det \phi_n}}\phi_n$.  This is a
section of $SL(E)$ over this annulus, satisfying
$$\dbar_{V_n}=\phi_n'^*\dbar_{A_n}.$$
It follows from the previous two lemmas that (after passing to a subsequence)
these gauge automorphisms converge to a limit $\phi'$ which satisfies
$$\dbar_V=\phi'^*\dbar_A.$$
On an annulus on which the winding number is odd, choose a holomorphic
function $\alpha$ from the annulus to $\cc^*$ with winding number 1, and let
$\psi=\left(\begin{array}{cc}
\alpha & 0\\
0 & 1 \end{array} \right).$
Then $\psi$ is a section of $GL(E)$ over $\Sigma\times W$.  The determinant of
$\psi\cdot\phi_n$ is even for all $n$, and
$$\dbar_{V_n}=(\psi\cdot\phi_n)^*(\psi^{-1})^*\dbar_{A_n}.$$
But then the previous two lemmas apply to the sequence of $SL$--automorphisms
$$\frac{1}{\sqrt{\det \psi\cdot\phi_n}}\psi\cdot\phi_n.$$
It follows that $\phi_n'=\frac{1}{\sqrt{\det \psi\cdot\phi_n}}\phi_n$
satisfy
$$\dbar_{V_n}=\phi_n'^*\dbar_{A_n}$$
on this annulus, and also that the sequence $\phi_n'$ has a convergent 
subsequence.

Summing up then, we see that we can find a sequence $\phi_n$
of sections of $GL(E)$ over $\Sigma\times W$ which satisfy
$$\dbar_{V_n}=\phi_n^*\dbar_{A_n}$$
on $\Sigma\times W$, and which have a subsequence which is converging in 
$C^\infty$ on compact subsets of $\Sigma\times W$ to a limit $\phi$ which
satisfies
$$\dbar_V=\phi^*\dbar_A.$$
The required operators $\dbar_n$ are then given by gluing $\dbar_{A_n}$
and $\dbar_{V_n}$ along $\Sigma\times W$ by the automorphisms $\phi_n$.

We may describe them in a way that shows clearly that they are extensions of
$\dbar_{A_n}$ on a fixed bundle ${\tilde E}$.  We know that
$\phi_n^{-1}\phi_N$ are converging to the identity in $C^\infty$ on
compact subsets of $\Sigma\times W$.
Choose $N$ large enough so that $\phi_n^{-1}\phi_N$ is homotopic to the
identity for all $n\ge N$.  Define ${\tilde E}$ to be the bundle formed
from $E\to\Sigma\times U$ and $E_V\to\Sigma\times V$ using the transition
function $\phi_N$ on $\Sigma\times W$.  Then
\begin{eqnarray*}
\dbar_n&=&\dbar_{A_n}\cup_{\phi_n}\dbar_{V_n}\\
       &=&\dbar_{A_n}\cup_{\phi_N}(\phi_n^{-1}\phi_N)^*\dbar_{V_n},
\end{eqnarray*}
for some extension of $\phi_n^{-1}\phi_N$ over $\Sigma\times V$, chosen so
that $\phi_n^{-1}\phi_N$ converge to the identity on compact subsets of
$\Sigma\times V$.
\endpf

We can now prove that
the two topologies on $\calm_{(\cald,c)}$ given in Section 
\ref{stablebundles} agree.

\begin{proposition}
\label{triplehomeo}
The correspondence given in Proposition \ref{triples} between \newline
$\calm_{(\cald,c)}$ and $\calm'_{(\cald,c)}$ is a homeomorphism.
That is to say, a sequence $[\cale_n]$ converges to $[\cale]$ in
$\calm_{(\cald,c)}$ if and only if the corresponding sequence of triples
\newline$(f_n,\cup_{i=1}^{k_n} z_{i;n},\cup_{i=1}^{k_n} [\cale_{i;n}])$
converges to the triple $(f_\cale,\cup_{i=1}^k z_i,\cup_{i=1}^k [\cale_i])$
associated to $[\cale]$ in the following way:
\begin{itemize}
\item $f_n \to f_\cale$ in $C^\infty$ on compact subsets of
$\cpone-\cup_{i=1}^k z_i$
\item for {\it any} open neighbourhood $D$ of $\cup_{i=1}^k z_i$, 
there exists $N$ such that \linebreak
$\cup_{i=1}^{k_n} z_{i;n}\subset D\quad \forall\quad n\ge N$
\item for {\it some} open neighbourhood $D$ of $\cup_{i=1}^k z_i$, the bundles
determined by $f_n$ and $\cup_{i=1}^{k_n} [\cale_{i;n}]$ converge in $C^\infty$
on $\Sigma\times D$ to that determined by $f_\cale$ and
$\cup_{i=1}^k [\cale_i]$.
\end{itemize}
 
\end{proposition}
\proof
Suppose first that
$[\cale_n]$ converges to $[\cale]$ in $\calm_{(\cald,c)}$.  We must show
that the associated sequence of triples satisfies the 3 conditions above.
The third requirement is satisfied immediately by taking $D=\cpone$.  Let us
consider the second.  Let $D$ be any open set in $\cpone$ containing
$\cup_{i=1}^k z_i$ and suppose that after passing to a subsequence at least
one of the points $\cup_{i=1}^{k_n} z_{i;n}$ is in the complement of $D$ for
all $n$.  Then a subsequence $z_{i_n;n}$ is converging to a limit point $z$ in
the compact set $\cpone-D$.  By continuity and the fact that stability is an
open condition, it follows that the limiting bundle $[\cale]$ must have
unstable restriction to $\Sigma\times \{z\}$, which is a contradiction.

It now follows that for any neighbourhood $D$ of $\cup_{i=1}^k z_i$ in $\cpone$
the restriction of $\cale_n$ to $\Sigma\times (\cpone-D)$ is stable on
every $\Sigma\times \{z\}$, for large $n$.  Thus the maps $f_n$ are converging
to $f_\cale$ in $C^\infty$ on $\cpone-D$, and hence in $C^\infty$ on
compact subsets of $\cpone-\cup_{i=1}^k z_i$.

Conversely, suppose
$(f_n,\cup_{i=1}^{k_n} z_{i;n},\cup_{i=1}^{k_n} [\cale_{i;n}])$
converges to 
$(f_\cale,\cup_{i=1}^k z_i,\cup_{i=1}^k [\cale_i])$.
Let $U$ be a disjoint union of open disks of radius $r$ centred at each of
$\cup_{i=1}^k z_i$, and let $U_\frac12\subset U$ be closed disks of radius 
$r/2$, also centred at $\cup_{i=1}^k z_i$.  Choose $r$ small enough so that
the restriction of $[\cale_n]$ is converging to the restriction of $[\cale]$
on $U$ (which we can do by the third condition in the definition of
convergence of triples).  By the second condition above, there exists
$N$ such that $\cup_{i=1}^{k_n} z_{i;n}\subset U_\frac12{}\!^\circ$ for all
$n\ge N$.

Let $V=\cpone-U_\frac12$.  Then $\cpone=U\cup V$, and $W=U\cap V$
is a union of open annuli.
We may choose $\dbar$--operators $\{\dbar_{A_n}\}_{n=1}^\infty,\dbar_A$
on $E\to\Sigma\times U$ so that
$$(E,\dbar_{A_n})\cong\cale_n, (E,\dbar_A)\cong\cale$$
on $\Sigma\times U$, and
$$\dbar_{A_n}\to\dbar_A$$
in $C^\infty$ on compact subsets of 
$\Sigma\times U$.

It follows from Corollary \ref{magic} that, after passing to a subsequence,
$$[\cale_n]\to[\cale].$$
In fact the same argument shows that any subsequence of $\{[\cale_n]\}$ has
a subsequence converging to $[\cale].$  Then since $\calm_{(\cald,\sigma)}$
has a  Hausdorff compactification, it follows that (for the entire sequence)
$$[\cale_n]\to[\cale],$$
as required.
\endpf

The following lemma is based on arguments used
by Donaldson \cite{d1} and Morgan \cite{m}, and in fact the proof closely
follows that in \cite[\S4.2]{m}.

\begin{lemma}
\label{donaldsonmorgan}
Let $X=\cup_{i=1}^l x_i$ be a finite set of points in $S$,
and let ${\tilde E}\to S$ be a smooth rank 2 hermitian bundle.
Suppose that $B_n$ are a sequence of unitary connections on ${\tilde E}$
with the following properties:
\begin{itemize}
\item $B_n$ converge to $B_\infty$ in $C^\infty$ on compact subsets of $S-X$,
where $B_\infty$ is a unitary connection on ${\tilde E}|_{S-X}$;
\item each $B_n$ is projectively ASD in some fixed neighbourhood $O$ of $X$ 
(with respect
to some fixed chosen \kahler metric on $S$).
\end{itemize}
Suppose each $({\tilde E},\dbar_{B_n})$ represents an element of
$\calm_{(\cald,c)}$, and $[({\tilde E},\dbar_{B_n})]$ converge to a limit
$[\xi]$ in $\calm_{(\cald,c)}$.  Suppose also that
$$({\tilde E}|_{S-X},\dbar_{B_\infty})\cong\zeta|_{S-X}$$
for some $\calc_\sigma$--stable bundle $\zeta\to S$.
Then $\xi\cong\zeta$.
\end{lemma}

\proof  The key point, as in \cite{d1} and \cite{m}, 
is to get uniform bounds on 
$|{\hat F}_{B_n}|$.  In both of the references just cited, these bounds
follow from the fact that the connections in question are ASD.  In fact though
since the sequence of connections are converging uniformly on compact subsets
of $S-X$, it is sufficient that they be ASD near $X$.

First change the connections $B_n$ and $B_\infty$ in the following way.
Give them all the same fixed central part.  Since they all determine
isomorphic holomorphic structures on the determinant line bundle 
this will not affect
the convergence on $S-X$ or the isomorphism classes of the holomorphic
structures they determine on ${\tilde E}$.

Since $[({\tilde E},\dbar_{B_n})]\to[\xi]$, it follows that there exist
$\dbar$--operators $\dbar_n'$,$\dbar_\infty'$ on ${\tilde E}$ with
${\rm lim}_{n\to\infty}\dbar_n'=\dbar_\infty'$, 
$({\tilde E},\dbar_\infty')\cong\xi$, and
\begin{equation}
\label{jn}
j_n\circ\dbar_{B_n}=\dbar_n'\circ j_n
\end{equation}
for some complex linear automorphism $j_n$ of ${\tilde E}$.
Let $B_n'$ denote the unique unitary connection on ${\tilde E}$ which is
compatible with the complex structure $\dbar_n'$ for each $n$.

\begin{claim} $|{\hat F}_{B_n}|$ and $|{\hat F}_{B_n'}|$
are bounded over $S$ uniformly in $n$.
\end{claim}
\proof
$|{\hat F}_{B_n'}|$ are uniformly bounded over $S$ since
$\dbar_n'$ are converging uniformly to $\dbar_\infty'$, and so 
$B_n'$ are also converging uniformly.

Similarly $|{\hat F}_{B_n}|$ are uniformly bounded on compact
subsets of $S-X$ since
$B_n$ is converging uniformly on $S-X$ to $B_\infty$.
$|{\hat F}_{B_n}|$ are uniformly bounded 
over $O$ since the $B_n$ are projectively ASD on $O$ with a fixed central
part.

Thus $|{\hat F}_{B_n}|$ are uniformly bounded over 
all of $S$.\endpf

Set $\tau_n=\tr(j_n^*\circ j_n)$, where $j_n^*$ is the adjoint of $j_n$
with respect to the metric on ${\tilde E}$.

\begin{claim} There exists a constant $C>0$ such that for all $n$
$${\rm sup}_{x\in S}\tau_n(x)
\le C\|\tau_n\|_{L^2(S)}.$$
\end{claim}
\proof This follows from the previous claim.  (See \cite[Corollary 4.2.9]{m}.)
\endpf

Multiply each $j_n$ by a constant $\lambda_n$ so that for all $n$
$$\|j_n\|_{L^4(S)}=\|\tau_n\|_{L^2(S)}=1.$$  
Then for all $n$ we have that
$${\rm sup}_{x\in S}\tau_n(x) \le C.$$
Now choose balls $B_i$ centred at each point $x_i\in X$ sufficiently small
so that they are disjoint from one another and so that the bundles $\xi$
and $\zeta$ are holomorphically trivial over each $B_i$.
Let $B=\cup_{i=1}^l B_i$ and let $T=S-B$.
The balls $B_i$ should also be chosen so that the volume of $B$ is at most
$\frac{1}{2C}$ (so that $\|j_n\|_{L^4(B)}\le\frac12$).

Write
\begin{eqnarray*}
\dbar_n'&=&\dbar_\infty'+\alpha_n\\
\dbar_{B_n}&=&\dbar_{B_\infty}+\beta_n
\end{eqnarray*}
where $\alpha_n,\beta_n$ converge to $0$ in $C^\infty$ on $T$.
Then (\ref{jn}) becomes
$$\dbar_\infty'\circ j_n-j_n\circ \dbar_{B_\infty}=j_n \beta_n-\alpha_n j_n.$$
Since the $\|j_n\|_{L^4(T)}$ are uniformly bounded for all $n$, it follows
from this equation and a standard bootstrap argument that there is a 
subsequence of the $j_n$ which converge in $C^\infty(T)$ to a limit $j_\infty$
which satisfies
$$\dbar_\infty'\circ j_\infty-j_\infty\circ \dbar_{B_\infty}=0.$$
That is to say $j_\infty$ defines a holomorphic map from $\zeta|_T$ to
$\xi|_T$.  Since $\|j_n\|_{L^4(T)}\ge\frac12$, it follows that
$j_\infty$ is not the zero map.  By making $B$ arbitrarily small,
$j_\infty$ can be extended to $S-X$.  Then by Hartogs' theorem $j_\infty$
extends to a nonzero holomorphic map from $\zeta$ to $\xi$.  Since these are
stable bundles with the same slope, this must in fact be an isomorphism.\endpf

The following result is based on Corollary \ref{magic} and 
Lemma \ref{donaldsonmorgan} and will be used to prove surjectivity of $\Psi$
and continuity of $\Psi^{-1}$.

\begin{lemma}
\label{nobubbles}
Let $\{[A_n]\}$ be gauge equivalence classes of smooth connections on 
\linebreak
\mbox{$E\to\cyl$} which are converging in $C^\infty$ on compact subsets on the 
complement of some finite set of points to some $[A]\in\calm$.  Suppose
that for any compact cylinder $U\subset S^1\times\rr$ there exists
$N_U$ such that $A_n|_{\Sigma\times U}$ is unitary and projectively
ASD for all $n>N_U$.
Let $\{[\xi_n]\},[\xi]$ be elements of $\calz_{(\cald,c)}$ with
$[\xi_n]\to[\xi]$ as $n\to\infty.$
Suppose also that
$$\xi_n|_{\cyl}\cong(E,\dbar_{A_n})$$
as holomorphic bundles.
Then $[\xi]=\Psi([A]).$
\end{lemma}

\proof 
Start by decomposing $\cpone$ into disjoint sets $U$ and $V$.  
Let $x_1,\dots,x_l$ be the ``bubble points'' in the
weak convergence of $[A_n]$ to $[A]$.
Choose $V$ to
be a sufficiently small neighbourhood of $\{0,\infty\}$ so that 
\begin{itemize}
\item for some $N$,
$\xi_n$ has stable restriction to $\Sigma\times\{z\}$ whenever $n\ge N$ and
$z\in V$.  (Such a $V$ must exist since $[\xi_n]$ are converging to a limit
in $\calz_{(\cald,c)}$.)
\item the points $x_1,\dots,x_l$ are all contained in $\Sigma\times(\cpone-V)$.
\end{itemize}
Let $U$ be a cylinder chosen so that $\cpone=U\cup V$ and $W=U\cap V$ is a
disjoint union of two annuli.

Let $f_n\co \cpone\to\caln(\Sigma)$ and $f_\infty\co \cpone\to\caln(\Sigma)$ be
the holomorphic maps associated to $\xi_n,\xi$ respectively
(as in Section \ref{stablebundles}).  Weak convergence
of $[A_n]$ to $[A]$ is sufficient to ensure that the holomorphic map
associated to $\Psi([A])$ is equal to $f_\infty$.

Note that gluing $(E,\dbar_{A_n})|_{\Sigma\times U}$ and
$({\rm Id}\times f_n)^*\calu|_{\Sigma\times V}$ along $W$ yields a bundle
isomorphic to $\xi_n$ for all $n$.  Also 
gluing $(E,\dbar_A)|_{\Sigma\times U}$ and
$({\rm Id}\times f_\infty)^*\calu|_{\Sigma\times V}$ along $W$ yields
a bundle in the equivalence class $\Psi([A])$.

Let ${\tilde E}\to S$ be the underlying smooth vector bundle of $\xi$ (and
hence also of $\xi_n$), and let $X=\cup_{i=1}^l x_i$.

Now by Corollary \ref{magic}, there exist $\dbar$--operators $\dbar_n$ 
on ${\tilde E}$ whose restriction to $U$ is equal to $\dbar_{A_n}$, and
a $\dbar$--operator $\dbar_\infty$ on ${\tilde E}|_{S-X}$
such that
\begin{itemize}
\item $({\tilde E},\dbar_n)\cong\xi_n$,
\item after passing to a subsequence, $\dbar_n\to\dbar_\infty$ in $C^\infty$ 
on compact subsets of $S-X$, and
\item $({\tilde E}|_{S-X},\dbar_\infty)\cong\Psi([A])|_{S-X}$.
\end{itemize}

Choose a metric $g$ on $S$ whose restriction to $U$ is equal to $g_{cyl}$, and
a hermitian metric $H$ on ${\tilde E}$ whose restriction to 
${\tilde E}|_{\Sigma\times U}$
is pulled back from $E_\Sigma\to\Sigma$.  Then we may form connections
$A(H,\dbar_n)$ on ${\tilde E}$, where $A(H,\dbar)$ 
denotes the unique
connection which is compatible with the metric $H$ and the complex structure
$\dbar$.

These connections satisfy the following properties:
\begin{itemize}
\item $A(H,\dbar_n)$ converge to $A(H,\dbar_\infty)$ in $C^\infty$ on
compact subsets of $S-X$
\item each $A(H,\dbar_n)$ is projectively ASD on $\Sigma\times U$ for all
$n\gg0$.
\end{itemize}

It follows from Lemma \ref{donaldsonmorgan} that $\xi\cong\Psi([A])$.\endpf

\section{From stable bundles to instantons}
\label{stabletoinstanton}
We will show that the map $\Psi$ is surjective by proving the
following result, where $e$ is some fixed element of $4\pi^2\zz$ and $\cald,c$
are given by (\ref{det}),(\ref{c2}) in Section \ref{defnot}.

\begin{proposition}
\label{psionto}
Let $\cale$ be a rank 2 holomorphic bundle on $S$ with 
$[\cale]\in\calz_{(\cald,c)}$.  Then there exists
a $U(2)$--instanton $A$ with $\Psi([A])=[\cale]$.
\end{proposition}

It may be possible to prove this using Donaldson's evolution equation method
directly, cf \cite{guo,d3}; however the degeneracy of the \kahler metric
$g_{cyl}$ on the 3--manifolds $\Sigma\times S^1 \times \{t\}$ may cause 
difficulties.  We avoid this analysis by finding a sequence of instanton
connections on $\cale\to S$ which converges to a $U(2)$--instanton on
$\Sigma \times S^1\times\rr$.  The proof uses two fundamental results: one
is Donaldson's theorem, which says that any bundle on an algebraic surface
which is stable with respect to a given Hodge metric admits a 
Hermitian--Einstein connection with respect to that Hodge metric.  
The other is Uhlenbeck's weak compactness for instantons.

Using these results it is quite straightforward to come up with a
$U(2)$--instanton $A$ on $E\to\Sigma \times S^1\times\rr$.  What is perhaps
surprising is that no energy is lost in the limiting process, and that
$\Psi([A])=[\cale]$.
\\

\noindent{\bf Proof of Proposition \ref{psionto}.}
We begin by finding a sequence of connections on $\cale$ which converge weakly
giving a $U(2)$--instanton $A$.

Consider the sequence of Hodge metrics on $S = \Sigma \times \cpone$ given by
$g_n = g_\Sigma + g_{\cpone,n}$, where
$$g_\Sigma = \; \mbox{chosen fixed \kahler metric on}\;  \Sigma, $$
$$g_{\cpone,1} = \; \mbox{Fubini--Study metric on}\;  \cpone. $$
In local coordinates $\theta, t$ on $S^1 \times \rr \subset \cpone$
$$g_{\cpone,1} = \frac{1}{4\pi\cosh^2t}(d\theta^2+dt^2).$$
We may choose real numbers $0=T_1<T_2<\dots,$ with $ T_n \rightarrow \infty$ as
$n \rightarrow \infty$, and smooth \kahler metrics $g_{\cpone,n}$ on 
$\cpone$
such that
$$g_{\cpone,n} = \left\{ \begin{array}{ll}
d\theta^2+dt^2 &
 \mbox{if} \quad t \in [-T_n+1,T_n-1]\\
\frac{1}{4\pi\cosh^2t}(d\theta^2+dt^2) &
 \mbox{if} \quad t \notin (-T_n,T_n)\\
\end{array} \right.$$
and the associated volume form $\omega_{\cpone,n}$ satisfies
$$\int_\cpone \omega_{\cpone,n} = n.$$
We will denote by $\omega_n \in \Omega^{1,1}(S,\rr)$ the \kahler form of the 
metric $g_n$.  Then $\omega_n$ is a representative of the cohomology class
$n\sigma+f$ in the \kahler cone $\calk(S)$ 
(see Section \ref{stablebundles}).  

The set $\calz_{(\cald,c)}$ is contained in $\calm_{(\cald,c)}$, so that
$\cale$ is stable with respect to any element of the chamber $\calc_\sigma$
in $\calk(S)$.
There is some $N \in \nn$, which depends on $e$, such that 
$[\omega_n] \in \calc_\sigma$ for all
$n \ge N$.

Let ${\tilde E}$ denote the underlying smooth vector bundle of $\cale$.  Let
$H_n$ be a hermitian metric on ${\tilde E}$ whose restriction to
${\tilde E}|_{\Sigma\times [-T_n+1,T_n-1]}$ is pulled back from 
$E_\Sigma\to\Sigma$.

Then by a theorem of Donaldson \cite {d1,dk}, 
for each $n \ge N$,
there is a unique unitary connection $A_n$ on $({\tilde E},H_n)$, up to gauge, 
which determines a holomorphic structure isomorphic
to $\cale$ and is
Hermitian--Einstein with respect to $g_n$.  This means that with respect to
$g_n$, $A_n$ is projectively anti-self-dual and $\tr F_{A_n}$ is a harmonic
representative for $-2\pi ic_1(\cale)$.

\begin{lemma}
The sequence  $\{[A_n]\}_{n \ge N}$ of connections on $\cale$ described above
has a subsequence which converges on compact subsets on the complement of a 
finite set of points in $\Sigma\times S^1\times\rr$ to give a (finite energy)
$U(2)$--instanton.
\end{lemma}
\proof The proof is a standard argument based on fundamental theorems of
Uhlenbeck.  See for example \cite[Theorem 6.1.1]{mmr}.\endpf

It now follows from Lemma \ref{nobubbles}, with $[\xi_n]=[\xi]=[\cale]$, that 
producing a $U(2)$--instanton in this way from the
holomorphic bundle $\cale$ actually produces an inverse of the map $\Psi$,
ie, that
$$\Psi([A])=[\cale].$$

\section{$\Psi$ is a homeomorphism}
\label{homeo}

In this section we complete the proof of Theorem \ref{mainthm} by showing
that the map
$$\Psi\co \calm_e/\cali\to\calz_{(\cald,c)}$$
is a homeomorphism.
We have already seen that it is a bijection; we need to see that both $\Psi$
and its inverse are continuous.

The codomain $\calz_{(\cald,c)}$ is first countable since it is 
a subspace of the Gieseker moduli space of $\calc_\sigma$--stable 
torsion-free sheaves.  
The domain is also first countable since it is contained in the Hilbert 
manifold of
finite-energy connections modulo gauge.  Thus to prove continuity it suffices
to show that $\Psi$ and $\Psi^{-1}$ preserve limits of sequences.

\begin{proposition}
\label{continuous}
$\Psi$ is continuous.
\end{proposition}
\proof
Suppose $[A_n]$ converge (strongly) to $[A]$ in $\calm_e$, for some 
$e\in 4\pi^2\zz$.
We must show that $\Psi([A_n])\to\Psi([A])$ in $\calm_{(\cald,c)}$.

Let $f_n\co \cpone\to\caln(\Sigma)$ and $f_\infty\co \cpone\to\caln(\Sigma)$ be
the holomorphic maps associated to $\Psi([A_n]),\Psi([A])$ respectively
(see Section \ref{stablebundles}).  Let $\cup_{i=1}^k z_i$ be the
marked points in $\cpone$ such that $\Psi([A])$ has unstable restriction
to $\Sigma\times \{z_i\}$.  Let $\cup_{i=1}^{k_n} z_{i;n}$ be the marked
points for $\Psi([A_n])$.

Let $D$ be any open neighbourhood of $\cup_{i=1}^k z_i$ in $\cpone$.  We will
show that there exists some $N$ such that the restriction
of $\Psi([A_n])$ to $\Sigma\times \{z\}$ is stable for all $z\notin D$ and for
all $n\ge N$.

For any compact set $K\subset S^1\times\rr$, it follows as in the proof of
Proposition \ref{triplehomeo} that there exists $N_K$ such that none of the 
marked points $\cup_{i=1}^{k_n} z_{i;n}$ are in $K-D$ for all $n\ge N_K$.
Thus either there exists an $N$ as required above, or some subsequence 
$z_{i_n;n}$ is converging to either $0$ or $\infty$.  But this would
contradict the strong convergence $[A_n]\to[A]$ since it would imply that
$\Psi([A])$ has a lower energy level than the limit of the
sequence $\Psi([A_n])$.

It also now follows (as in the proof of Proposition \ref{triplehomeo})
that $f_n\to f_\infty$ in $C^\infty$ on $\cpone-\cup_{i=1}^k z_i$.

By Proposition \ref{triplehomeo}, it now follows that
$\Psi([A_n])\to\Psi([A])$ in $\calm_{(\cald,c)}$ as required.\endpf

\begin{proposition}
\label{continuousinverse}
$\Psi^{-1}$ is continuous.
\end{proposition}
\proof
Suppose $\{[\xi_n]\},[\xi]$ are elements of $\calz_{(\cald,c)}$ with
$[\xi_n]\to[\xi]$ as $n\to\infty$.  By Propositions \ref{injective} and
\ref{psionto} there exist unique elements $[A_n]\in\calm_e/\cali$ with
$\Psi([A_n])=[\xi_n]$.  Uhlenbeck compactness implies that
 the sequence $[A_n]$ has a subsequence
which is converging on compact sets on the complement of
a finite set of points ${x_1,\dots,x_l}$ to $[A]$, where $A$ is a 
$U(2)$--instanton with $e(A)\le e-4\pi^2l$.  By Lemma \ref{nobubbles},
$$[\xi]=\Psi([A]).$$
It follows that in fact $[A]\in\calm_e/\cali$.
The same argument shows that any subsequence of ${[A_n]}$ has a subsequence
converging to $[A]$; it follows that $[A_n]$ converge to $[A]$ as required.
\endpf

\begin{remark} One could follow the argument in \cite[\S4]{m} and show that 
if in
fact bubbling occurs in the weak convergence $[A_n]\to[A]$, then the limit
of the sequence $\Psi([A_n])$ will have sheaf 
singularities at the bubble points.
This partially describes the extension of $\Psi$ (or rather $\Psi^{-1}$) to the
compactifications of $\calm_e$ and $\calz_{(\cald,c)}$.
\end{remark}

\begin{remark}
It is to be expected that the map $\Psi$ may carry more information.  By
considering some suitable weighted Sobolev completions it should be possible
to compare the deformation complexes and Kuranishi models associated to
$[A]$ and $\Psi([A])$, where $[A]\in\calm_e$, following 
\cite[Sections 4.3.3, 4.3.4]{fm}, and hence show that $\Psi$ is an 
orientation-preserving isomorphism of real analytic ringed spaces.
\end{remark}

\end{document}